\documentclass[10pt]{article} 

\usepackage{latexsym,amsmath,amssymb,a4wide,amscd}
\usepackage[latin1]{inputenc}
\usepackage[T1]{fontenc}

\usepackage{array}

\usepackage{amsthm}

\usepackage{url}

\usepackage{mathabx}

\usepackage{MnSymbol}

\usepackage{longtable}

\numberwithin{equation}{section}
\numberwithin{figure}{section}

\swapnumbers
\theoremstyle{definition}
\newtheorem {satz}{Theorem}[section]
\newtheorem {conj}[satz]{Conjecture}
\newtheorem {remark}[satz]{Remark}
\newtheorem {example}[satz]{Example}
\newtheorem {deff}[satz]{Definition}
\newtheorem {lemma}[satz]{Lemma}
\newtheorem {prop}[satz]{Proposition}
\newtheorem {kor}[satz]{Corollary}
\newtheorem {question}[satz]{Question}

\usepackage[pdftex]{graphicx} 

\usepackage[nooneline]{caption}

\begin{document}

\date{}

\title{\Large {\bf Normal surfaces as combinatorial slicings}}

\author{Jonathan Spreer}

\maketitle

\subsection*{\centering Abstract}

{\em
We investigate slicings of combinatorial manifolds as properly embedded co-dimension $1$ submanifolds. A focus is given to dimension $3$ where slicings are (discrete) normal surfaces. For the cases of $2$-neighborly $3$-manifolds as well as quadrangulated slicings, lower bounds on the number of quadrilaterals of slicings depending on its genus $g$ are presented. These are shown to be sharp for infinitely many values of $g$. Furthermore we classify slicings of combinatorial $3$-manifolds which are weakly neighborly polyhedral maps.
}\\
\\
\textbf{MSC 2000: } {\bf 57Q15};  
57M20 \\ 
\\
\textbf{Keywords: } normal surface, slicing, combinatorial manifold, weakly neighborly polyhedral map, Heegaard genus, combinatorial Heegaard splitting.

\section{Discrete normal surfaces and slicings}
\label{sec:normsurf}

\begin{figure}[htp]
	\centering
	\includegraphics[width=0.9\textwidth]{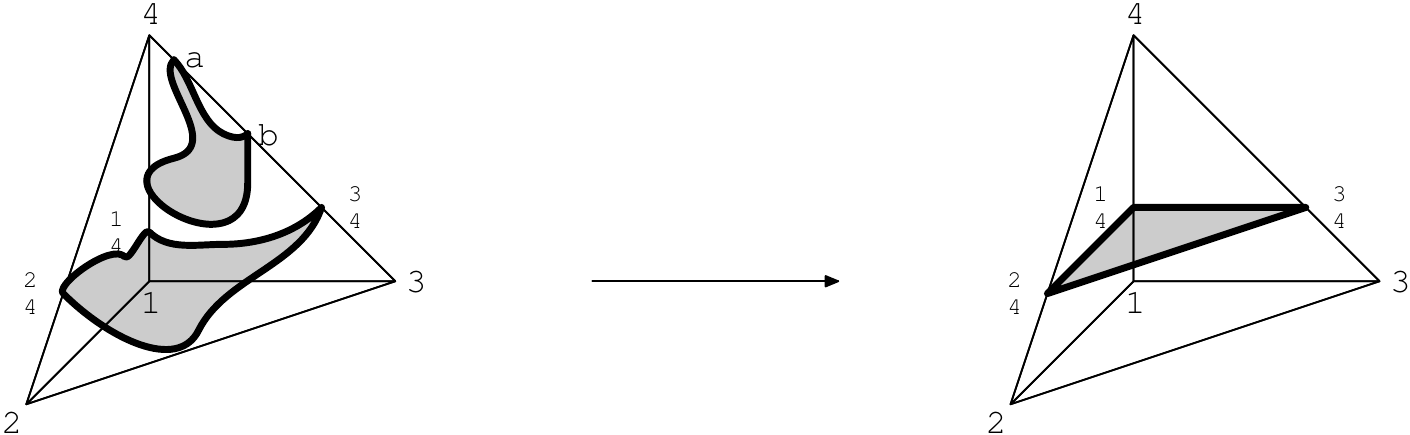}
	\caption{Intersection of an embedded surface with a tetrahedron of the surrounding combinatorial $3$-manifold and the corresponding normal subset.}
	\label{fig:kneser}
\end{figure}

In this article we develop a combinatorial theory of discrete normal surfaces in combinatorial $3$-manifolds. The concept of {\it normal surfaces} is due to Kneser \cite{Kneser29ClosedSurfIn3Mflds}. He used it to prove one part of the prime decomposition theorem in the theory of $3$-manifolds. A surface $ S $, properly embedded into a $3$-manifold $M$, is said to be {\it normal}, if it respects a given cell decomposition of $M$ in the following sense: It does not intersect any vertex nor touch any 3-cell of the manifold and does not intersect with any 2-cell in a circle or an arc starting and ending in a point of the same edge (see Figure \ref{fig:kneser} for the simplicial case).

The precise definition of the term {\it normal surface} is due to Haken \cite{Haken61TheoNormFl}. Haken developed an algebraic theory of normal surfaces to advance the research on the homeomorphism problem of $3$-manifolds (for any pair of $3$-manifolds $(M_1,M_2)$ decide in a finite number of steps whether $M_1 \cong M_2 $ or not, cf. \cite{Haken62HomeomProb3Mflds}). In the theory of (hyperbolic) $3$-manifolds, normal surfaces are often examined using special kinds of cell decompositions: If $\tilde{\Delta}$ is a set of tetrahedra together with a set of gluing instructions $\Phi$ on the set of triangles of $\tilde{\Delta}$ such that each triangle is identified with at most one other triangle, then $P = \tilde{\Delta} / \Phi$ is called a {\it pseudo triangulation}. 

However, in this article we consider only {\it combinatorial manifolds}: A combinatorial $d$-manifold ({\it combinatorial $d$-pseudomanifold}) $M$ is a $d$-dimensional, pure, simplicial complex whose vertex links are all combinatorial spheres with standard PL-structure (combinatorial manifolds). Note that every combinatorial manifold (pseudomanifold) is also a pseudo triangulation. The $f${\it -vector} of a combinatorial manifold (pseudomanifold) $M$ is a $(d+1)$-tuple of integers $f(M)$ where the $i$-th entry $f_{i-1}$ denotes the number of $(i-1)$-dimensional faces of $M$. We call $M$ $k${\it -neighborly}, if $f_{k-1} = { f_0 \choose k }$, i. e. if $M$ contains all possible ($k-1$)-dimensional faces.

\begin{deff}[Polyhedral manifold, polyhedral map, cf. \cite{Brehm93PolyhedralMflds}]
	A {\em polyhedral complex} $C$ is a finite family of convex polytopes such that {\it (i)} for every polytope $P \in C$ all of its faces $F \in P$ are contained in $C$ and {\it (ii)} the intersection $P_1 \cap P_2$ of any two polytopes $P_1, P_2 \in C$ is either empty or a common face of $P_1$ and $P_2$.
	
	A {\em polyhedral manifold} is a polyhedral complex $M$ such that there exists a simplicial subdivision of $M$ which is a combinatorial manifold. If $M$ is a surface we will call it a {\em polyhedral map}. If, in addition $M$ entirely consists of $m$-gons, we call it a {\it polyhedral $m$-gon map}.
\end{deff}

\begin{figure}[ht]
	\centering
	\includegraphics[width=\textwidth]{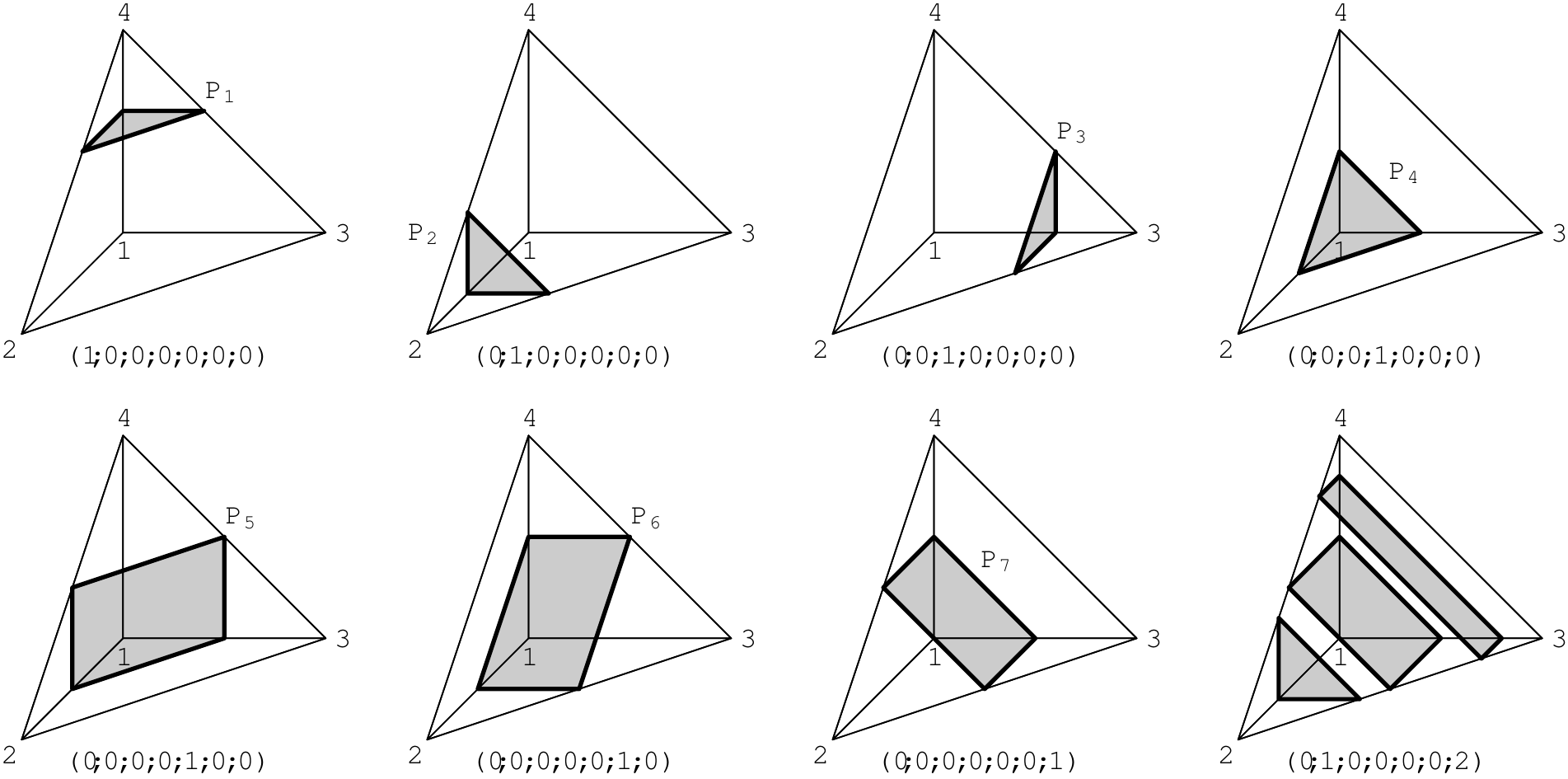}
	\caption{Possible cuts of a tetrahedron by a normal surface and associated normal coordinates. Note that the leftmost picture of the bottom row can not be part of a discrete normal surface.}
  \label{fig:haken}
\end{figure}

\begin{deff}[Discrete normal surface]
	
	\noindent
	Let $M$ be a combinatorial $3$-manifold ($3$-pseudomanifold), $\Delta \in M$ one of its tetrahedra and $P$ the intersection of $\Delta$ with a plane that does not include any vertex of $\Delta$. Then $P$ is called a {\it normal subset} of $\Delta$. Up to an isotopy that respects the face lattice of $\Delta$, $P$ is equal to one of the triangles $P_{i}$, $1 \leq i \leq 4$, or quadrilaterals $P_{i}$, $5 \leq i \leq 7$, shown in Figure \ref{fig:haken}. 
	
	A polyhedral map $S \subset M$ that entirely consists of facets $P_{i}$ such that every tetrahedron contains at most one facet is called {\it discrete normal surface} of $M$.
\end{deff}

\begin{remark}
	In classical normal surface theory a tetrahedron of a combinatorial manifold can contain several facets of possibly different types. Hence, we can assign a vector $v \in \mathbb{N}_0^7$ to each tetrahedron counting the number of parallel cuts of each type. The set of all such vectors of all tetrahedra of a combinatorial manifold is called the {\it normal coordinates set} of a normal surface (cf. Figure \ref{fig:haken}). Since all cuts of a tetrahedron have to be disjoint, not all combinations of cuts are valid. In particular, each tetrahedron may contain only one type of quadrilateral at most. In addition, the cuts of any two adjacent tetrahedra have to be compatible as a normal surface has to be a closed polyhedral map. For a combinatorial $3$-manifold $M$ we have $3 f_2 $ such linear {\it compatibility equations} where $f_2$ denotes the number of triangles of $M$ ($3$ restrictions per triangle). A normal coordinates set with compatible entries is called {\it admissible}. 	It is an interesting fact that a normal surface without vertex linking connected components is already determined by its quadrilaterals (see \cite{Tillmann07NormalSurfacesTopFinite3Mflds}, Thm. 2.4). This leads to a more compressed type of normal coordinates with a vector $v \in \mathbb{N}_0^3$ for each tetrahedron.

	The description of normal surfaces in terms of normal coordinates gives rise to the concept of the {\em geometric sum}: The union of two normal surfaces, defined by the componentwise sum $\Sigma$ of their normal coordinates, is well defined if and only if $\Sigma$ is admissible (see \cite{Tillmann07NormalSurfacesTopFinite3Mflds} for further details). This implies that the theory of normal surfaces has a lot less algebraic structure than homology theory, although Haken himself emphasized a close connection between these two theories in \cite{Haken61TheoNormFl}. 
\end{remark}

\bigskip
Let us now introduce the notion of \textit{slicings} of combinatorial manifolds:
It is well known from classical Morse theory, that the pre-image of a non-critical point of a smooth Morse function on a closed smooth $3$-manifold is a properly embedded closed surface. In the field of $PL$-topology K\"uhnel developed what one might call a polyhedral Morse theory  (compare \cite{Kuehnel90TrigMnfFewVert,Kuehnel95TightPolySubm}):

\begin{deff}[Rsl-function \cite{Kuehnel95TightPolySubm}]
 \label{def:rsl}
 Let $M$ be a combinatorial $d$-manifold. A function $f:M \to \mathbb{R}$ is called \textit{regular simplexwise linear (rsl)} if $f(v) \neq f(w)$ for any two vertices $w \neq v$ of $M$ and if $f$ is linear when restricted to an arbitrary simplex of $M$.
 
 A point $x \in M$ is said to be {\it critical} for a rsl-function $f:M \to \mathbb{R}$ if \[H_{\star} (M_x , M_x \backslash \{ x \} , F) \neq 0 \] where $M_x := \{ y \in M \, | \, f(y) \leq f(x) \}$ and $F$ is a field. Here $H_{\star}$ denotes an appropriate homology theory.
\end{deff}

It follows that no point of $M$ can be critical except possibly the vertices.

\begin{deff}[Slicing]
	
	\label{def:slicing}
	\noindent
	Let $M$ be a combinatorial pseudomanifold of dimension $d$ and $f:M \to \mathbb{R}$ a rsl-function. Then we call the pre-image $f^{-1} (x)$ a {\it slicing} of $M$ whenever $x \neq f(v)$ for any vertex $v \in M$.
\end{deff}

By construction, a slicing is a polyhedral $(d-1)$-manifold and for any ordered pair $x < y$ we have $f^{-1} (x) \cong f^{-1} (y)$ whenever $f^{-1}([x,y])$ contains no critical vertex of $M$. In particular, a slicing $S$ of a closed combinatorial $3$-manifold $M$ is a discrete normal surface: It follows from the simplexwise linearity of $f$ that the intersection of the pre-image with any tetrahedron of $M$ either forms a single triangle or a single quadrilateral. In addition, if two facets of $S$ lie in adjacent tetrahedra they either are disjoint or glued together along the intersection line of the pre-image and the common triangle.

\begin{remark}
	Any partition $V = V_1 \dot{\cup} V_2 $ of the set of vertices of $M$ already determines a slicing: Just define a rsl-function $f: M \to \mathbb{R}$ with $f(v) < f(w)$ for all $v \in V_1$ and $w \in V_2$ and look at a pre-image $f^{-1} (x_0)$ for any $f(v) < x_0 < f(w)$. In the following we will write $S_{(V_1,V_2)} := f^{-1} (x_0)$ for the slicing defined by the vertex partition $V = V_1 \dot{\cup} V_2 $.
	
	Every vertex of a slicing is given as an intersection point of the corresponding pre-image with an edge $\langle u,w \rangle$ of the combinatorial manifold. Since there is at most one such intersection point per edge, we usually label this vertex of the slicing according to the vertices  of the corresponding edge, that is $\binom{u}{w}$ with $u \in V_1$ and $w \in V_2$.
	
	By construction, every slicing decomposes the surrounding combinatorial manifold $M$ into at least $2$ pieces (an upper part $M^+$ and a lower part $M^-$). This is not the case for discrete normal surfaces in general. However, in the sequel we will focus on discrete normal surfaces that are slicings and we will apply the above notation for discrete normal surfaces whenever this is possible.
\end{remark}

\begin{figure}[ht]
	\includegraphics[width=\textwidth]{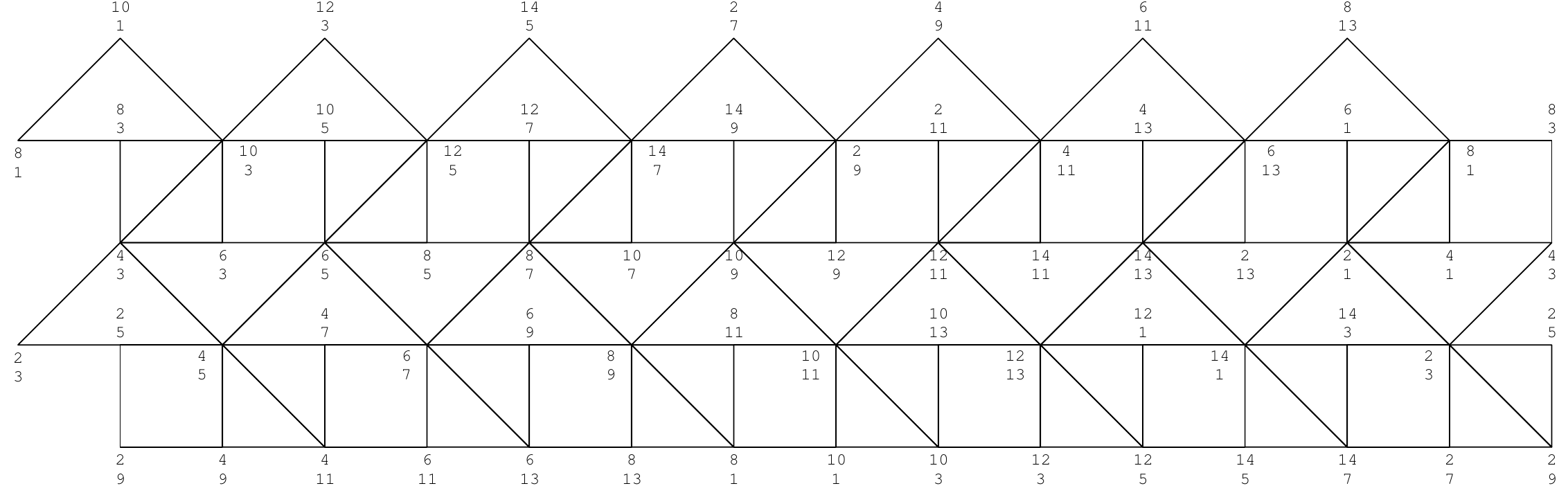}
	\caption{Slicing of genus $1$ with 35 Quadrilaterals, 28 triangles, $7 \cdot 7 = 49$ vertices and cyclic $\mathbb{Z}_{14}$-symmetry. It it is obtained from a $2$-neighborly $14$-vertex version of the lens space $\operatorname{L}(3,1)$ with transitive automorphism group.}
	\label{fig:lense}	
\end{figure}

Since every combinatorial pseudomanifold $M$ has a finite number of vertices, there exist only a finite number of slicings of any fixed $M$. Hence, if $f$ is chosen carefully, the induced slicings admit a useful visualization of $M$. This has been done already in a number of publications: See \cite{Kuehnel84RhombiTess3Space15Vertex3Torus} for a visualization of a $15$-vertex version of the $3$-torus, \cite{Spreer09CombPorpsOfK3} for some $3$-dimensional slicings of the Casella-K\"{u}hnel triangulation of the $K3$-surface and \cite{Kuehnel95TightPolySubm} for various further examples.  Figure \ref{fig:lense} shows the separating torus of a $2$-neighborly $14$-vertex triangulation of the lens space $\operatorname{L}(3,1)$ with transitive automorphism group (triangulation $^3 14^{1}_{6}$ in \cite{Lutz03TrigMnfFewVertVertTrans}).

\bigskip
\noindent
The paper is organized as follows:

In Section \ref{sec:Topology} we investigate upon the possible topological types of discrete normal surfaces a given combinatorial $3$-manifold admits. In particular, we present a minimal combinatorial Heegaard splitting of the $3$-torus and small spheres separating two arbitrary bounded $3$-manifolds.

In Section \ref{sec:strucSli} we discuss the local combinatorial structure of slicings of combinatorial $3$-manifolds ($3$-pseudomanifolds). A variety of observations on the different roles of triangles and quadrilaterals is presented. 

In Section \ref{sec:trigquad} we present our main result which is a lower bound on the number of quadrilaterals of a slicing depending on its genus and assuming certain properties such as $2$-neighborliness of the surrounding manifold. Furthermore we discuss some cases of equality. 

In Section \ref{sec:neighslic} we examine slicings which are weakly neighborly polyhedral maps. A condition for the weakly neighborliness is given, as well as a classification of all weakly neighborly slicings of combinatorial manifolds. 

\section{The genus of discrete normal surfaces}
\label{sec:Topology}

For embedded orientable surfaces $S \subset M$ which decompose a $3$-manifold $M$ into two pieces we have the following statement:

\begin{prop}
	\label{thm:mvseq}
	Let $M$ be a connected compact orientable $3$-manifold and $S \subset M$ a properly embedded connected surface decomposing $M$ into two bounded $3$-manifolds $M^-$ and $M^+$ with common boundary $S$. Then
	\begin{equation}
		\label{eq:mvseq}
		\beta_1 ( M^- ; \mathbb{Z} ) - \beta_2 ( M^- ; \mathbb{Z} ) = g ( S ) =	\beta_1 ( M^+ ; \mathbb{Z} ) - \beta_2 ( M^+ ; \mathbb{Z} ),
	\end{equation}
	holds, where $g(S)$ denotes the genus of $S$ and $\beta_i ( M^{\pm} ; \mathbb{Z} )$ the $i$-th integral Betti number of $M^{\pm}$.
\end{prop}

\begin{proof} 
Since $M^+$ ($M^-$) is a bounded manifold with orientable boundary $S$ (the orientability of $S$ follows from the orientability of $M$ and the fact, that $S$ is separating  $M^+$ from $M^-$) we can glue $M^+$ ($M^-$) with a copy of itself along its boundary $S$ obtaining a closed $3$-manifold $M$. Recall that $\chi (M) = 0$ since $M$ is a closed $3$-manifold. By the additivity of the Euler characteristic it now follows that
\begin{equation}
	0 = \chi (M) = 2 \chi (M^+) - \chi (S) = 2(1 - \beta_1 (M^+) + \beta_2 (M^+)) - ((2 - \beta_1 (S))
\end{equation}
which directly leads to 
\begin{equation}	
	 \beta_1 (M^+) - \beta_2 (M^+) \, = \, \frac12 \beta_1 (S) \, = \, g ( S ).
\end{equation}
The calculation for $M^-$ is the same.
\end{proof}

Since the embedding of $S$ is arbitrary, the genus $g(S)$ does not depend on any properties of $M$. In particular any $3$-manifold $M$ admits an embedding of a connected orientable surface $S$ of any genus $g(S)$. 

Of course, if $M$ is a combinatorial manifold and $S$ is a discrete normal surface the situation is somewhat different. Due to the finite number of tetrahedra the genus of an embedded discrete normal surface $S$ is always bounded. In fact the only topological type of discrete normal surfaces that occurs in any combinatorial manifold $M$ is the $2$-sphere (for example as the vertex figure of an arbitrary vertex). In the sequel we will investigate restrictions on the genus of $S$ given by the combinatorial properties of $M$.

\begin{prop}
	\label{prop:handles}
	Let $M$ be an orientable connected combinatorial $3$-manifold, $V = V_1 \dot{\cup} V_2$ a partition of the set of vertices and $|V| \in \{ 2n, \, 2n+1 \}$, $n \in \mathbb{N}$, such that $S_{(V_1,V_2)}$ is connected. Then
	\begin{equation}
		\label{eq:handles}
		g(S_{(V_1,V_2)}) \leq { n-1 \choose 2 }.
	\end{equation}
\end{prop}

\begin{proof}
	Let $M = M^- \cup M^+$, $V_1 \subset M^+$, $V_2 \subset M^-$, be the decomposition of $M$ with common boundary $S_{(V_1,V_2)}$. From Proposition \ref{thm:mvseq} it follows that 
	\begin{equation*}
		\beta_1 (M^+) - \beta_2 (M^+) = g(S_{(V_1,V_2)}).
	\end{equation*}
	Now let $e$ be the number of edges, $t$ the number of triangles, $\Delta$ the number of tetrahedra of $\operatorname{span}(V_1)$ and w.l.o.g. $m:= |V_1| \leq |V_2|$. Then
	\begin{eqnarray}
		\label{eq:eulerhandles}
		g(S_{(V_1,V_2)}) &=& 1 - 1 + \beta_1 (M^+) - \beta_2 (M^+) \nonumber \\
	  &=& 1 - m + e - t + \Delta
	\end{eqnarray}
	holds by the Euler-Poincar\'e formula. Since $t \geq 2 \Delta$, $e \leq { m \choose 2 }$ and $m \leq n$, the right hand side of (\ref{eq:eulerhandles}) is maximal for $t = \Delta = 0$, $e = { m \choose 2 }$ and $m = n$, and thus $ g(S_{(V_1,V_2)}) \leq 1 - n + { n \choose 2} = { n-1 \choose 2 }$.
\end{proof}

See Figure \ref{fig:slicing4-5} (right) for an example attaining equality in \ref{eq:handles} in the case $n=5$ and $g(S_{(V_1,V_2)})=6$. 

In order to make a closer connection between the genus of $S$ and the topology of $M$, we will restrict ourselves to special kinds of decompositions in the following:

\begin{enumerate}
	\item A {\em Heegaard splitting} \cite[Definition 12.11]{Lickorish} decomposes a manifold $M$ into two {\em handle bodies} $M^+$ and $M^-$, i. e. two $3$-balls with a number of solid cylinders attached in a handle-like manner with an orientable surface $S$ as common boundary. The number of handles is equal for both $M^+$ and $M^-$ and is called their {\em genus} coinciding with the genus of $S$ (a fact that also follows from Proposition \ref{thm:mvseq}). It is also referred to as the {\em genus of the Heegaard splitting}. Every $3$-manifold admits such a splitting \cite[Lemma 12.12]{Lickorish}.
	\item A {\em Heegaard splitting of minimal genus} is a Heegaard splitting of $M$ where the resulting handle bodies do not have more than the minimum number of handles needed. Such a splitting defines the {\em Heegaard genus $g(M)$}, a topological invariant of $M$. Thus, a slicing defining such a minimal Heegaard splitting partially reflects the topological structure of $M$.
	\item A slicing $S_{(V_1,V_2)}$ induced by a vertex splitting $V = V_1 \dot{\cup} V_2 $ such that the underlying set of $\operatorname{span} (V_i)$, $i = 1,2$, defines bounded $3$-manifolds is called {\it separating surface}. Examples of such discrete normal surfaces do not put any restrictions on the topology of $M$ in both $M^+$ and $M^-$. In particular $M$ can be extended to a closed combinatorial $3$-manifold $\hat{M}$ in which $S_{(V_1,V_2)}$ is no longer separating $\hat{M}$ into two pieces $\hat{M}^{\pm}$. This gives rise to a family of examples of discrete normal surfaces that possibly are no slicings (in $\hat{M}$).
\end{enumerate}

\begin{figure}[ht]
	\centering
	\includegraphics[width=0.5\textwidth]{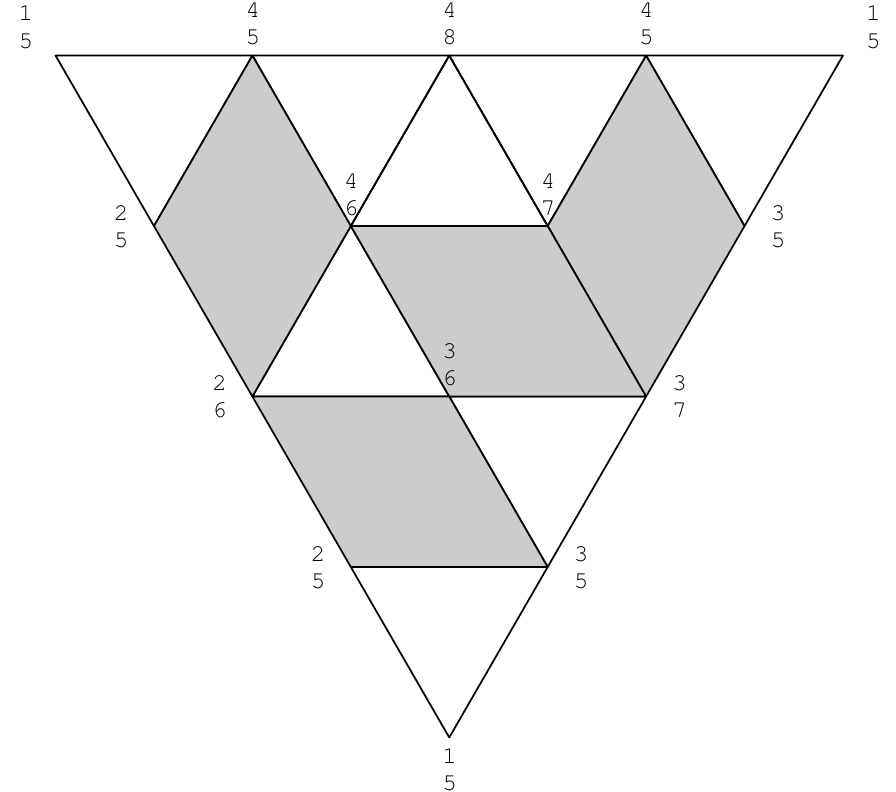}
	\caption{Slicing of $C_1$ in between $\{ 1,2,3,4 \}$ and $ \{ 5,6,7,8 \}$ with the minimum number of facets.}
	\label{fig:minslice}
\end{figure}

\begin{figure}[ht]
	\centering
	\includegraphics[width=0.5\textwidth]{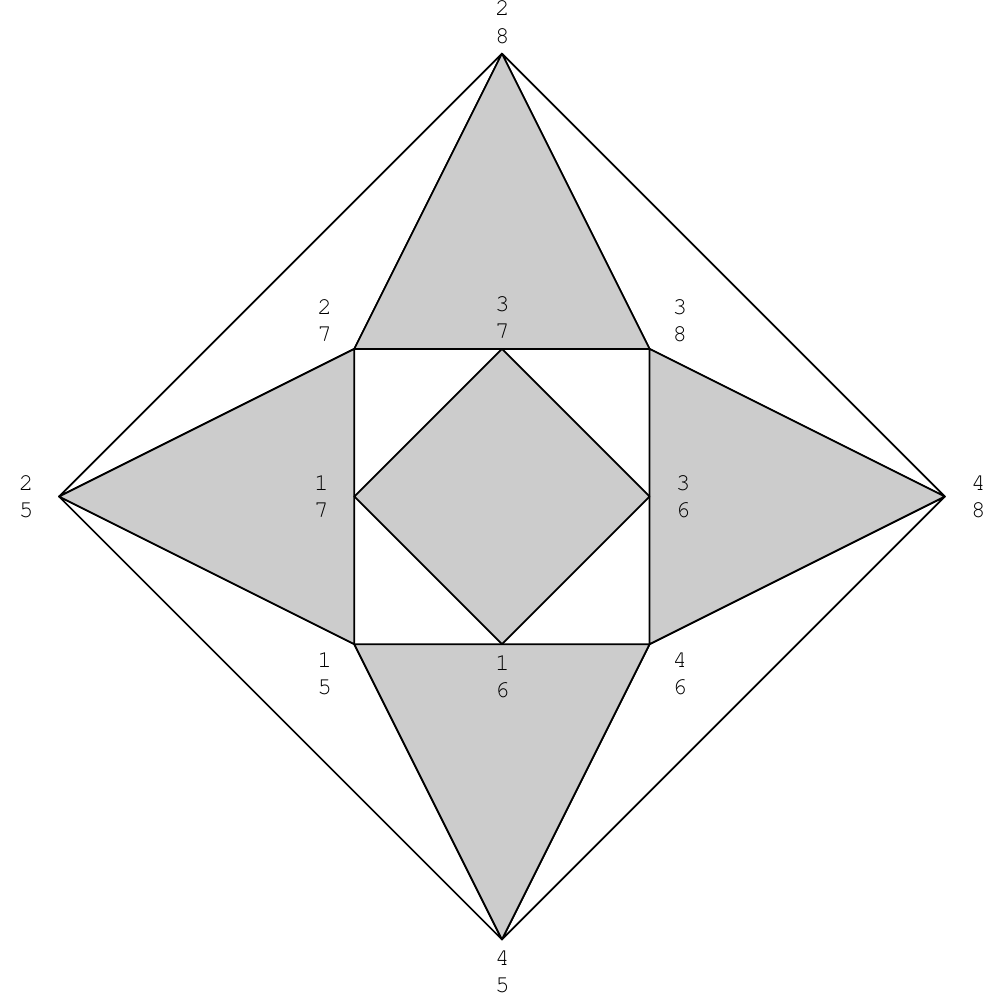}
	\caption{Schlegel diagram of the slicing of $C_2$ in between $\{ 1,2,3,4 \}$ and $\{ 5,6,7,8 \}$ (a cuboctahedron).}
	\label{fig:cuboctahedron}
\end{figure}

\begin{example}[Separating surfaces]
	Consider a combinatorial version of the cylinder $C = S_g \times [0,1]$ where $S_g$ is a triangulated oriented surface of genus $g$. C is a bounded $3$-manifold that can be closed by adding an arbitrary number of combinatorial cells on both sides. A slicing through $C$ disjoint to its boundary can be seen as a separating surface of genus $g$. 
	
	In the case $g=0$ we can consider the cylinder
	\begin{center}
		\begin{tabular}{l@{}l@{}l@{}l@{}l@{}l@{}l}
			$ C_1 = \langle $&
			$\langle 1\,2\,3\,5 \rangle$, &
			$\langle 1\,2\,4\,5 \rangle$, &
			$\langle 1\,3\,4\,5 \rangle$, &
			$\langle 2\,3\,4\,6 \rangle$, &
			$\langle 2\,3\,5\,6 \rangle$, &
			$\langle 2\,4\,5\,6 \rangle$,\\
			&
			$\langle 3\,4\,5\,7 \rangle$, &
			$\langle 3\,4\,6\,7 \rangle$, &
			$\langle 3\,5\,6\,7 \rangle$, &
			$\langle 4\,5\,6\,8 \rangle$, &
			$\langle 4\,5\,7\,8 \rangle$, &
			$\langle 4\,6\,7\,8 \rangle \rangle $
		\end{tabular}
	\end{center}
	\noindent
	with boundary 
	\begin{center}
		\begin{tabular}{l@{}l@{}l@{}l@{}l}
			$ \partial C_1 = \langle $ &
			$\langle 1\,2\,3 \rangle$, &
			$\langle 1\,2\,4 \rangle$, &
			$\langle 1\,3\,4 \rangle$, &
			$\langle 2\,3\,4 \rangle$,\\
			&
			$\langle 5\,6\,7 \rangle$, &
			$\langle 5\,6\,8 \rangle$, &
			$\langle 5\,7\,8 \rangle$, &
			$\langle 6\,7\,8 \rangle \rangle $.
		\end{tabular}
	\end{center}
	Since we need at least $8$ vertices for a triangulation of the boundary of a cylinder of type $C = S^2 \times [0,1]$, we need at least $8$ vertices to triangulate $C$. Barnette's Lower Bound Theorem (cf. \cite{Barnette71MinNumVertSimplePoly,Barnette73ProofLBCConvPoly}) tells us that  $C$ needs at least $3 \cdot 8-10-2 = 12$ tetrahedra. Now let us consider a slicing $S$ through $C$ with $\partial(C) \cap S = \emptyset $. If there exist a set of tetrahedra $\Delta \subset C$ disjoint to $S$, we define $\tilde{C} = C \setminus \Delta$.  By the arguments above $\tilde{C}$ again has at least $12$ tetrahedra and all of them intersect with $S$. Thus, $S$ must have at least $12$ facets, $8$ of then being triangles. Figure \ref{fig:minslice} shows a slicing of $C_1$ with $8$ triangles and the minimum number of $4$ quadrilaterals. The slicing is a subdivided tetrahedron. 
	
	A slightly different cylinder with $8$ vertices 
	\begin{center}
		\begin{tabular}{l@{}l@{}l@{}l@{}l@{}l@{}l@{}l}
			$C_2 = \langle $&
			$\langle 1\,2\,3\,7 \rangle$, &
			$\langle 1\,2\,4\,5 \rangle$, &
			$\langle 1\,2\,5\,7 \rangle$, &
			$\langle 1\,3\,4\,6 \rangle$, &
			$\langle 1\,3\,6\,7 \rangle$, &
			$\langle 1\,4\,5\,6 \rangle$, &
			$\langle 1\,5\,6\,7 \rangle$, \\
			&
			$\langle 2\,3\,4\,8 \rangle$, &
			$\langle 2\,3\,7\,8 \rangle$, &
			$\langle 2\,4\,5\,8 \rangle$, &
			$\langle 2\,5\,7\,8 \rangle$, &
			$\langle 3\,4\,6\,8 \rangle$, &
			$\langle 3\,6\,7\,8 \rangle$, &
			$\langle 4\,5\,6\,8 \rangle \rangle $
		\end{tabular}
	\end{center}
	\noindent
	and necessarily the same boundary as $C_1$ leads to a more symmetric slicing in form of a cuboctahedron with $14$ facets, necessarily $8$ triangles and $6$ quadrilaterals (see Figure \ref{fig:cuboctahedron}).
	
	$C_1$ was obtained by canonically subdividing the prism-complex $S_4^2 \times [0,1]$ into $3 \cdot 4 = 12 $ tetrahedra. This procedure is available for any closed oriented surface $S_g$ and gives rise to a slicing through $S_g \times [0,1]$ with $f_2(S_g)$ quadrilaterals and $2 \cdot f_2(S_g)$ triangles. The question whether or not this type of subdivision contains the minimum number of quadrilaterals needed for arbitrary values of $g$ in the case that $S_g$ is a vertex minimal triangulation is interesting but does not seem to be answered yet. It is closely related to a generalization of Theorem \ref{thm:quadrsli}. An equivalent reformulation would be:
	
	\begin{question}
		Let $C:= S_g \times [0,1]$ be a triangulation of a thickened orientable surface of genus $g$ such that no vertex lies in the interior of $C$. What is the minimum number of tetrahedra that do not share a triangle with the boundary of $C$?
	\end{question}
\end{example}

\begin{example}[Combinatorial Heegaard splittings]
	It is well known that the three-torus $\mathbb{T}^3$ is a $3$-manifold of Heegaard genus three. By a {\em combinatorial Heegaard splitting} we mean the decomposition of a combinatorial $3$-manifold into two polyhedral handle bodies such that the common boundary is a slicing. It seems likely that the minimal genus of a combinatorial Heegaard splitting equals the Heegaard genus of the underlying manifold, although this does not seem to be trivial. Hence, a visualization of a minimal combinatorial splitting of the three-torus would be interesting.
	
	\begin{figure}
		\centering
		\includegraphics[width=0.5\textwidth]{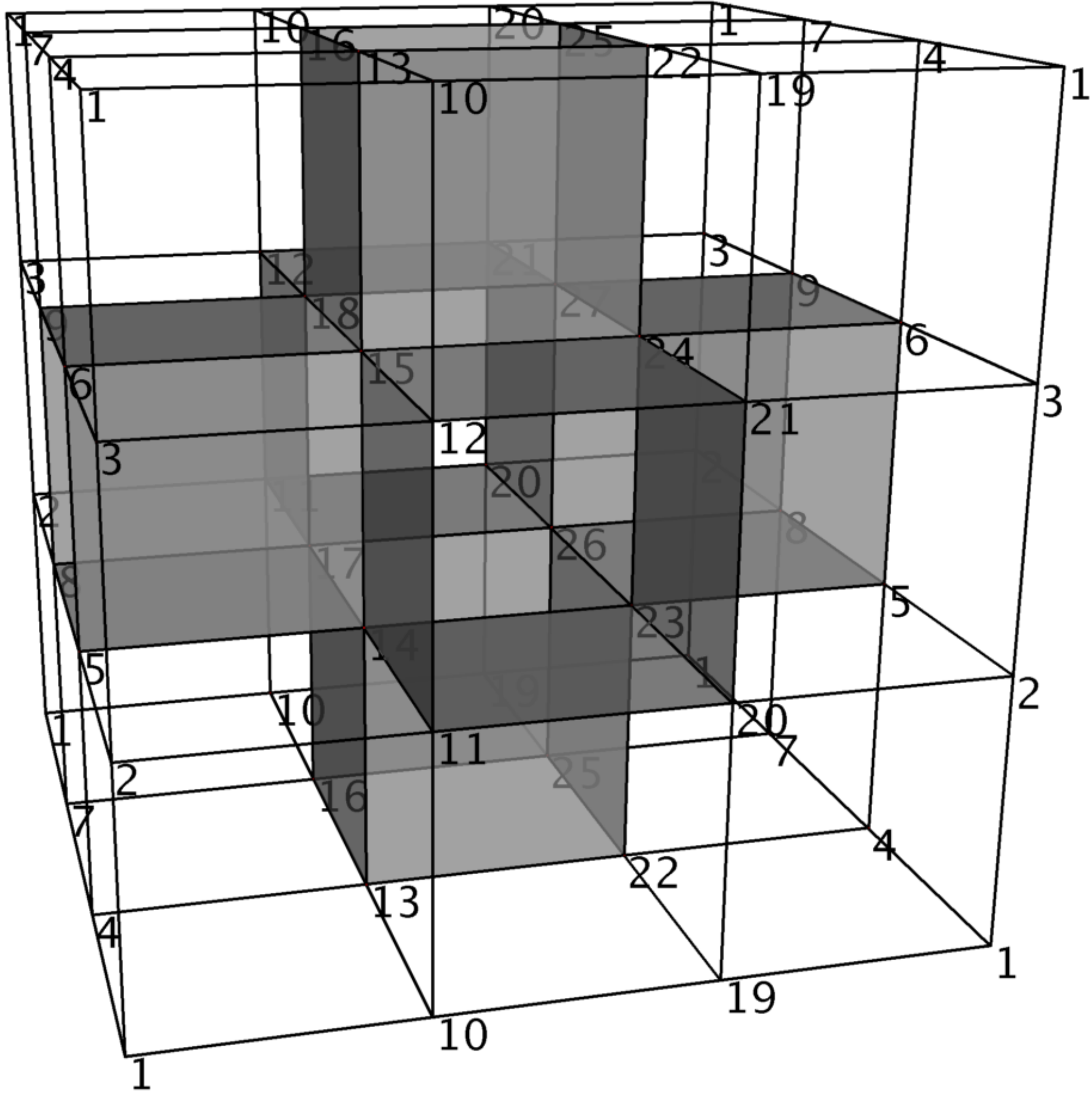}
		\caption{Combinatorial Heegaard splitting of the canonical $3 \times 3 \times 3$-cube subdivision of the $3$-torus. Opposite faces are identified.}
		\label{fig:heegaard3}
	\end{figure}
	
	The decomposition of $\mathbb{T}^3$ as a $3 \times 3 \times 3$ subdivided cube with a pairwise identification of opposite faces and $3^3 = 27$ vertices admits a splitting of (minimal) genus $3$ as shown in Figure \ref{fig:heegaard3}. We can easily transform the cube decomposition into a combinatorial manifold by subdividing each cube into $6$ tetrahedra. The embedded surface of genus three is a sub-complex and, thus, not a slicing. However, by a slight perturbation such an example is attained. Note, that the perturbed slicing would be too complex to admit a useful visualization in this context. 
\end{example}

\section{The combinatorial structure of slicings}
\label{sec:strucSli}

As already stated, it is difficult to link global properties of a discrete normal surface or, in higher dimensions, a slicing to its surrounding combinatorial manifold (pseudomanifold). However, we will observe that its local combinatorial structure in fact depends on the manifold. 

\begin{remark}
	Slicings are very special polyhedral manifolds. In order to see this let $P$ be a facet of a slicing $S$ of a combinatorial $d$-manifold ($d$-pseudomanifold) $M$. Then  
	\begin{equation}
		\label{eq:facets}
		P \cong_{\operatorname{comb.}} \Delta^{d-1-k} \times \Delta^k ,
	\end{equation}
	$ 0 \leq k \leq \frac{d-1}{2}$ and $S$ contains at most $\lfloor \frac{d+1}{2} \rfloor$ different types of polytopes as facets. In particular, (\ref{eq:facets}) implies that a $2$-dimensional slicing (or a discrete normal surface) only consists of triangles $\Delta^2$ and quadrilaterals $\Delta^1 \times \Delta^1$. Hence, the search for relations between the number of triangles and quadrilaterals in the $2$-dimensional case seems natural and has been investigated by Kalelkar in \cite{Kalelkar08ChiQuadrNS}. We will return to that question in Section \ref{sec:trigquad}.	
\end{remark}

If the slicing of a combinatorial $3$-manifold (pseudomanifold) $M$ is a vertex figure of $M$, it obviously contains no quadrilaterals. Hence, any triangulated sphere (surface) can be seen as the vertex figure of a suitable combinatorial $3$-manifold ($3$-pseudomanifold). However, every slicing different from a disjoint union of vertex figures contains quadrilaterals as facets, since in this case both $M^-$ and $M^+$ at least contain one edge. Thus, every connected slicing of a combinatorial $3$-manifold with genus $>0$ has to contain quadrilaterals whereas this is not true for combinatorial $3$-pseudomanifolds.

\begin{figure}[ht]
	\centering
	\includegraphics[width=0.5\textwidth]{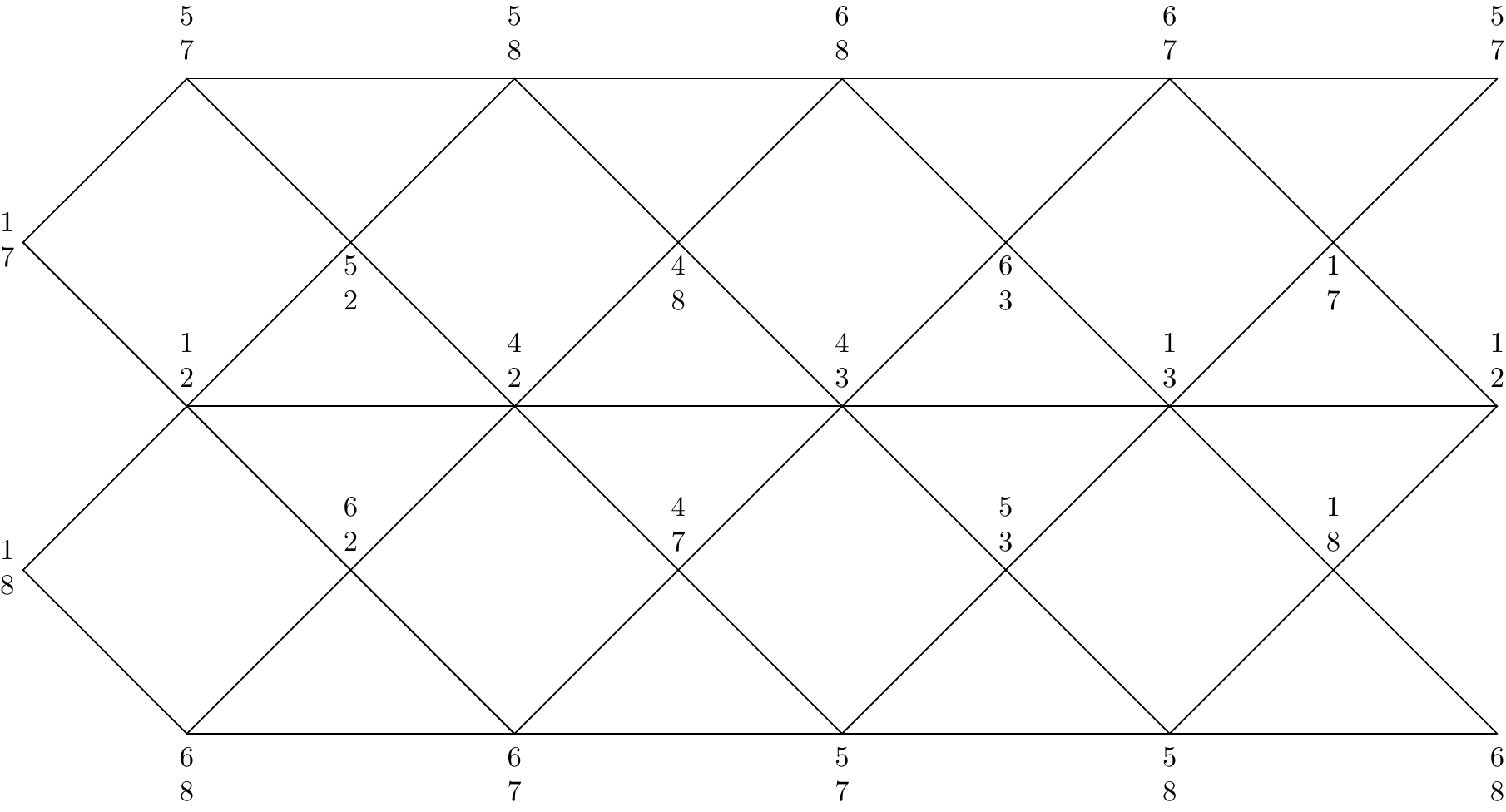}
	\caption{Highly symmetric centered slicing of genus $1$ through an $8$-vertex triangulation of a $3$-dimensional Kummer variety.}
	\label{fig:kummer}
\end{figure}

It seems to be fact that not only the vertex figures of singular vertices of combinatorial pseudomanifolds are different from slicings of combinatorial manifolds: For example, consider the slicing $S_{((1,4,5,6),(2,3,7,8))} \cong \mathbb{T}^2$ through the $8$-vertex triangulation of a $3$-dimensional (singular) Kummer variety $K^3$ (from \cite{Kuhnel86MinTrigKummVar}) shown in Figure \ref{fig:kummer}. The $24$ tetrahedra of the triangulation are completely determined by $S_{((1,4,5,6),(2,3,7,8))}$.	Note, that there exists a basis $ \langle \alpha , \beta \rangle $ of $ H_1 (S_{((1,4,5,6),(2,3,7,8))})$ for which  both cycles are quadrilaterals of the form	$ \langle \langle { a \choose c }, { a \choose d } \rangle , \langle { a \choose d }, {b \choose d } \rangle , \langle {b \choose d },  {b \choose c } \rangle , \langle {b \choose c } , { a \choose c } \rangle \rangle $. If we look at $\alpha$ and $\beta$ in $M^{+}$ they both collapse to an edge and are thus both contractible. The same holds for $\alpha$ and $\beta$ in $M^{-}$. In fact, we have $\operatorname{span} (V_i) \cong S^2_4$ for $i=1,2$. This contradicts with Proposition \ref{thm:mvseq}. As a consequence,	$S_{((1,4,5,6),(2,3,7,8))}$ can not be a slicing of a combinatorial manifold.

Despite these differences the following definitions and lemma hold for both combinatorial $3$-manifolds and $3$-pseudomanifolds. In the following, let $S_{(V_1,V_2)}$ be the slicing defined by the vertex splitting $V = V_1 \dot{\cup} V_2$ of a combinatorial $3$-pseudomanifold.

\begin{deff}
	\label{def:trace}
	Let $S_{(V_1,V_2)} \subset M$ be a slicing of a combinatorial $3$-pseudomanifold $M$ defined by  $V = V_1 \dot{\cup} V_2$ and $x \in V_1$ a vertex. We define
	\begin{eqnarray}
		C^x_2 &:=& \{ \delta \in S_{(V_1,V_2)} \, | \, \delta = \left\langle \textstyle \binom{x}{a}, \binom{x}{b}, \binom{x}{c} \right\rangle; \,  a,b,c \in V_2 \} \\
		C^x_1 &:=& \{ \delta \in S_{(V_1,V_2)} \backslash C^x_2 \, | \, \delta = \left\langle \textstyle \binom{x}{a}, \binom{x}{b} \right\rangle; \, a,b \in V_2 \} \\
		C^x_0 &:=& \{ \delta \in S_{(V_1,V_2)} \backslash (C^x_2 \cup C^x_1) \, | \,  \delta = \left\langle \textstyle \binom{x}{a} \right\rangle; \, a \in V_2 \}
	\end{eqnarray}
	and
	\begin{equation*}
		C^x := \overline{C^x_2 \cup C^x_1}.
	\end{equation*}
	We call $C^x$ the \emph{trace of} $x$ \emph{on} $S_{(V_1,V_2)}$. Analogously we define $C_y$ for any $y \in V_2$.
\end{deff}

Note, that $ C^x \cup C^x_0 = S_{(V_1,V_2)} \cap S_{(\{x\},V \backslash \{x\})} $ where $S_{(\{x\},V \backslash \{x\})} $ is the vertex figure of $x$ which uniquely consists of triangles.

\begin{lemma}
	\label{lemma:neighbors}
	Let $S_{(V_1,V_2)} \subset M$ be a slicing of a combinatorial $3$-pseudomanifold $M$, $x,y \in V_1$, $a,b \in V_2$. If $\left\langle \binom{x}{a} , \binom{y}{b} \right\rangle$ is an edge of $ S_{(V_1,V_2)} $, then either $a=b$ or $x=y$ holds.
\end{lemma}

\begin{proof}
	If $\left\langle \binom{x}{a}, \binom{y}{b} \right\rangle$ is an edge of a slicing $ S_{(V_1,V_2)} $ in a $3$-manifold $M$ it must be a subset of some triangle $ \langle \alpha, \beta, \gamma \rangle \in M$. So $\binom{x}{a}$ and  $\binom{y}{b}$ are both center points of edges in $\langle \alpha, \beta, \gamma \rangle$. It follows immediately that $a$, $b$, $x$ and $y$ can not all be different.
\end{proof}

\begin{kor}
	\label{cor:form}
	All quadrilaterals $Q$ of a slicing $ S_{(V_1,V_2)} $ of a combinatorial $3$-pseudomanifold $M$ are of the form
	\begin{equation*}
		\textstyle Q = \left\langle \binom{x}{a}, \binom{x}{b}, \binom{y}{a}, \binom{y}{b} 
		\right\rangle
	\end{equation*}
	and all triangles $T$ are of the forms
	\begin{equation*}
		\textstyle T = \left\langle  \binom{x}{a}, \binom{x}{b}, \binom{x}{c} \right\rangle
		\quad \textrm{ or } \quad T = \left\langle 
		\binom{x}{a}, \binom{y}{a}, \binom{z}{a}
		\right\rangle
	\end{equation*}
	with $x,y,z \in V_1$, $a,b,c \in V_2 $.
\end{kor}

\begin{proof}
	Recall that all facets of $ S_{(V_1,V_2)} $ originate from proper sections with a tetrahedron.
\end{proof}

\begin{lemma}
	\label{lemma:intersection}
	Let $S_{(V_1,V_2)} \subset M$ be a slicing of a combinatorial $3$-pseudomanifold $M$ defined by  $V = V_1 \dot{\cup} V_2$ and $C^x$ ($C_a$) the trace of $x \in V_1$ ($a \in V_2$). Then the following implications hold:
	\begin{eqnarray}
	i) & x,y \in V_1, x \neq y & \Rightarrow \quad C^x \cap C^y = \emptyset \\
	ii) & x \in V_1, a \in V_2 & \Rightarrow \quad C^x \cap C_a = \left\lbrace \begin{array}{l}
	\left\langle \binom{x}{a} \right\rangle 
	\textrm{ if } \langle x,a \rangle \in M\\
	\emptyset \textrm{ otherwise} \end{array} \right.
	\end{eqnarray}
\end{lemma}

\begin{proof} 
	This follows immediately from the definition of $C^x$:
	
	\medskip
	Since $x,y \in V_1$, the upper entries of all vertices in $C^x$ ($C^y$) are equal  to $x$ ($y$). But $x \neq y $ and the intersection of $ C^x $ and $C^y$ must be empty. This proves $i)$.
	
	\medskip
	With the same argument we see that $C^x \cap C_a \subseteq \left\langle \binom{x}{a} \right\rangle$. However, the vertex $\binom{x}{a}$ is part of $S_{(V_1,V_2)}$ iff $\langle x,a \rangle$ is an edge of $M$. This shows $ii)$.
\end{proof}

\begin{lemma}
	\label{lemma:intersection2}
	Let $Q$ be a quadrilateral of a slicing $S_{(V_1,V_2)}$ of a combinatorial $3$-pseudomanifold $M$, $C^x$ the trace of a vertex $x\in V_1$ in $S_{(V_1,V_2)}$ and $C_a$ the trace of $a\in V_2$. Then $Q$ shares at most one edge with $C^x$ and $C_a$, respectively.
\end{lemma}

\begin{proof} 
	By Corollary \ref{cor:form}, any quadrilateral $Q$ of $S_{(V_1,V_2)}$ is of the form $ Q = \left\langle \binom{x}{a}, \binom{x}{b}, \binom{y}{a}, \binom{y}{b} \right\rangle $. This implies that $Q$ shares exactly one edge with $C^x$, $C^y$, $C_a$ and $C_b$.
\end{proof}

As already mentioned, discrete normal surfaces without vertex linking components are completely determined by their quadrilaterals. The observations about the local combinatorial structures of slicings made above emphasize this crucial property.  

\section{Triangles vs. quadrilaterals}
\label{sec:trigquad}

Discrete normal surfaces are polyhedral maps consisting of triangles and quadrilaterals. Although triangulated discrete normal surfaces exist, we already pointed out that any discrete normal surface of a combinatorial $3$-manifold with non-trivial genus has to contain quadrilateral facets and thus is not simplicial. The relation between the genus of a normal surface and the number of quadrilaterals was investigated by Kalelkar in \cite{Kalelkar08ChiQuadrNS}:

\begin{satz}[Kalelkar, \cite{Kalelkar08ChiQuadrNS}]
	\label{thm:kalelkar1}
	Let $S$ be a closed, oriented, connected normal surface of a pseudo triangulation of a closed $3$-manifold $M$, $g$ its genus and $q$ the number of quadrilaterals in $S$, then we have
	\begin{equation} \label{eq:kalelkar1} g \leq \frac72 q. \end{equation}
\end{satz}

Although Theorem \ref{thm:kalelkar1} holds for arbitrary values of $g$ it merely assures the appearance of one quadrilateral for normal surfaces of genus $\leq 3$ and two quadrilaterals in the case $3 < g \leq 7$. But in practice we have seen that usually a lot more quadrilaterals occur. Hence, (\ref{eq:kalelkar1}) does not seem to be sharp in the case of discrete normal surfaces and combinatorial $3$-manifolds. In Chapter \ref{sec:strucSli} we investigated some rules describing how strongly connected components of triangles determine areas of quadrilaterals in the case of slicings. Moreover, we know from Proposition \ref{prop:handles} that an increasing genus of a slicing yields an increasing complexity of the manifold, thus an increasing minimum number of vertices and, in turn, an increasing complexity of the slicing. In order to improve inequality (\ref{eq:kalelkar1}) in a combinatorial setting we will start with some universal observations on combinatorial $3$-manifolds before imposing some restrictions on the generic case:

Let $M$ be an orientable combinatorial $3$-manifold with $f$-vector $f$ and $V$ the set of vertices. From the Lower Bound Theorem (LBT, cf. \cite{Kalai87RigidityLBT}) for combinatorial manifolds we get the following restrictions on the number of edges $f_1$ (here the rightmost inequality is just the trivial upper bound $ f_1 \leq \binom{f_0}{2}$): 
\begin{equation}
	\label{eq:lbt}
	4f_0 - 10 \leq f_1 \leq \binom{f_0}{2}.
\end{equation}
In addition,the Dehn-Sommerville equations for $3$-manifolds hold:
\begin{eqnarray}
	f_0 - f_1 + f_2 - f_3 & = & 0 \nonumber \\
	2 f_2 - 4 f_3 & = & 0. \nonumber
\end{eqnarray}
Now, let $S_{(V_1,V_2)}$ be a slicing defined by the vertex partition $ V = V_1 \dot{\cup} V_2 $ and let $n$ be the number of vertices of $S_{(V_1,V_2)}$. Then the obvious condition 
\begin{equation}
	\label{eq:2neighb}
	n \leq |V_1||V_2|
\end{equation}
holds, with equality whenever every vertex of $V_1$ is connected to every vertex of $V_2$ by an edge. Moreover, we have the following equalities which we call {\it Dehn-Sommerville equations for slicings}:
\begin{eqnarray}
  n - e + t + q & = & 2 - 2g, \nonumber \\
  - 2 e + 3 t + 4 q & = & 0, \nonumber
\end{eqnarray}
where $e$ is the number of edges, $t$ the number of triangles, $q$ the number of quadrilaterals and $g$ the genus of $S_{(V_1,V_2)}$. Note that in the following we will call the $4$-tuple $(v,e,t,q)$ the $f$-vector of a slicing despite the fact, that the $f$-vector of $S_{(V_1,V_2)}$ seen as a polyhedral map would be the $3$-tuple $(v,e,t+q)$.

Finally, we have another obvious relation between $M$ and $S_{(V_1,V_2)}$: Since every tetrahedron of $M$ contains at most one facet of $S_{(V_1,V_2)}$ it is
\begin{equation}
   \label{eq:link}
   f_3 - t - q \geq 0,
  \end{equation}
and more precisely the number $f_3 - t - q \geq 0$ equals the number of tetrahedra in the span of $V_1$ and $V_2$. For a fixed $n$, these equations induce a linear system of equations of dimension $4 \times 7$ with rank $4$:
\begin{center}
	\begin{tabular}{ccccccc|l}
		$f_1$ & $f_2$ & $f_3$ & $ v $ & $ e $ & $ t $ & $ q $ & \\	\hline
		$-1$     & $1$      & $-1   $  &       &       &       &       & $ -n $ \\
				     & $-2$     & $4    $  &       &       &       &       & $ 0 $ \\
					   &          &          &$ 1 $  & $ -1 $& $ 1 $ & $ 1 $ & $ 2 - 2g$\\
					   &          &          &       & $ -2 $& $ 3 $ & $ 4 $ & $ 0 $. \\
	\end{tabular}
\end{center} 
Solutions of these equations additionally have to fulfill the inequalities (\ref{eq:lbt}), (\ref{eq:2neighb}) and (\ref{eq:link}).

Note that we can not expect useful information from the above linear system in the general case (a fact that is also implied by the separating surfaces from Section \ref{sec:Topology}). However, if we restrict $M$ to be $2$-neighborly, the situation is different.

In the following, let $S_{(V_1,V_2)}$ be a slicing of a $2$-neighborly combinatorial $3$-manifold $M$ such that $|V_1| + c = |V_2| - c$ for a $c \in \frac12 \mathbb{N}$. In this case we get a $4 \times 4$ system of the form
\begin{center}
	\begin{tabular}{cccc|l}
		\label{tab:comb}
		$ f_3 $ & $ e $ & $ t $ & $ q $ & \\
		\hline
		$ 1 $      &       &       &       & $ \binom{n}{2} - n $ \\
						   & $ -1 $& $ 1 $ & $ 1 $ & $ 2 - 2g - \frac{n^2}{4} + c^2 $ \\
						   & $ -2 $& $ 3 $ & $ 4 $ & $ 0 $ \\
		$ 1 $      &       & $ -1$ & $ -1$ & $ \geq 0 $. \\
	\end{tabular}
\end{center} 
As the system is of rank $4$, $n$ and $c$ completely determine the combinatorial properties of the discrete normal surface up to the relation $f_3 \geq t + q$.

From this we can deduce our main result:

\begin{satz}
	\label{thm:result}
	Let $S_{(V_1,V_2)}$ be a slicing of a $2$-neighborly combinatorial $3$-manifold $M$ and let $g$, $n$, $q$ and $c$ be defined as above. Then we have
	\begin{equation}
		q \geq 4g + \frac{3n}{2} - (4 + 2c^2).
		\label{eq:result} 
	\end{equation}
	For $c=0$, Inequality (\ref{eq:result}) is sharp for all values $g = \binom{\frac{n}{2}-1}{2}$.
\end{satz}

\begin{proof}
	The row echelon form from the linear system above is 
	\begin{center}
		\begin{tabular}{cccc|l}
			 $ f_3 $ & $ e $ & $ t $ & $ q $ & \\ \hline
			$ 1 $ &     &       &       & $ \binom{n}{2} - n $ \\
						&$-1$ & $ 1 $ & $ 1 $ & $ 2 - 2g - \frac{n^2}{4} + c^2 $ \\
						&     & $ 1 $ & $ 2 $ & $ 4g - 4 + \frac{n^2}{2} - 2c^2 $ \\
						&     &       & $ 1 $ & $ \geq 4g + \frac{3n}{2} - (4 + 2c^2) $. \\
		\end{tabular}
	\end{center} 
	\noindent
	This proves (\ref{eq:result}).
	
	Now let $c=0$. This implies that $n = 2k$ for some $k \in \mathbb{Z}$ and thus $g = \binom{k-1}{2}$ is well defined. By replacing the value of $g$ with $\binom{k-1}{2}$ in (\ref{eq:result}) we get
	\begin{equation*}
		q \geq 4\binom{k-1}{2} + 3k - 4 = {2k \choose 2} - 2k = f_3 . 
	\end{equation*}
  The last step follows from the Dehn-Sommerville equations for combinatorial $3$-manifolds and (\ref{eq:link}), which implies that $q = 4g + 3k - 4$.

	In order to show the tightness of (\ref{eq:result}) we construct a series of slicings with the above properties: Consider the ($2$-neighborly) boundary complex $BdC_4(2k)$ of the cyclic $4$-polytope $C_4(2k)$ (cf. \cite{Gruenbaum67ConvPoly,Kuehnel85NeighbComb3MfldsDihedralAutGroup}).
	
	For any $k$ the boundary complex $BdC_4(2k)$ can be constructed as follows: Given the dihedral group in the following permutation representation
	\begin{equation*}
		D_{2k} = \left \langle (1, \ldots , 2k ) , (2k,2)(2k-1,3) \ldots (k +2,k) \right \rangle
	\end{equation*}
	with $4k$ elements, take the union of the following $(k-2)$ orbits of length $2k$
	\begin{equation*}
		(1,2,3,4)_{2k} , (1,2,4,5)_{2k} , \ldots , (1,2,k,k +1)_{2k} 
	\end{equation*}
	and the orbit $(1,2,k+1,k+2)_{k}$ of length $k$. The resulting complex $C_{2k}$ is a $2$-neighborly combinatorial $3$-sphere with the maximum number of $ 2k(k-2) + k = { 2k \choose 2 }- 2k$ facets.
	
	As one easily deduces from the group action, all of the facets of $C_{2k}$ contain exactly two even and exactly two odd vertex labels. This property is also known as {\em Gale's evenness condition} of cyclic polytopes (see for example \cite{Gruenbaum67ConvPoly} page $62$). If one defines $V_1 = \{ 1,3,\ldots ,2k-1 \}$ and	$V_2 = \{ 2,4,\ldots ,2k \}$, neither $\operatorname{span} (V_1)$  nor $\operatorname{span} (V_2)$ contains a triangle, and the induced slicing is an example of equality in (\ref{eq:result}).
\end{proof}

\begin{example}
	The simplest case for equality in (\ref{eq:result}) is the equilibrium set in the boundary complex of the join of two triangles which, in turn, is nothing else than
	\begin{center}
		\begin{tabular}{l@{}lll}
			$ BdC_4(6) = \langle $ &
			$\langle 1\,2\,3\,4 \rangle$, &
			$\langle 1\,2\,3\,6 \rangle$, &
			$\langle 1\,2\,4\,5 \rangle$,\\
			&
			$\langle 1\,2\,5\,6 \rangle$, &
			$\langle 1\,3\,4\,6 \rangle$, &
			$\langle 1\,4\,5\,6 \rangle$,\\
			&
			$\langle 2\,3\,4\,5 \rangle$, &
			$\langle 2\,3\,5\,6 \rangle$, &
			$\langle 3\,4\,5\,6 \rangle \rangle$.
		\end{tabular}
	\end{center}
	\noindent
	By separating the odd from the even vertex labels as above, one obtains the slicing
	\begin{center}
		\begin{tabular}{l@{}lll}
			$ S_{(\{ 1,3,5 \}, \{ 2,4,6 \})} = \langle $ &
			$\langle {1 \choose 2}, {1 \choose 4}, {3 \choose 2}, {3 \choose 4} \rangle$, &
			$\langle {1 \choose 2}, {1 \choose 6}, {3 \choose 2}, {3 \choose 6} \rangle$, &
			$\langle {1 \choose 2}, {1 \choose 4}, {5 \choose 2}, {5 \choose 4} \rangle$, \\[0.2cm]
			&
			$\langle {1 \choose 2}, {1 \choose 6}, {5 \choose 2}, {5 \choose 6} \rangle$, &
			$\langle {1 \choose 4}, {1 \choose 6}, {3 \choose 4}, {3 \choose 6} \rangle$, &
			$\langle {1 \choose 4}, {1 \choose 6}, {5 \choose 4}, {5 \choose 6} \rangle$, \\[0.2cm]
			&
			$\langle {3 \choose 2}, {3 \choose 4}, {5 \choose 2}, {5 \choose 4} \rangle$, &
			$\langle {3 \choose 2}, {3 \choose 6}, {5 \choose 2}, {5 \choose 6} \rangle$, &
			$\langle {3 \choose 4}, {3 \choose 6}, {5 \choose 4}, {5 \choose 6} \rangle \rangle ,$
		\end{tabular}
	\end{center}
	\noindent	
	which is the standard $3 \times 3$-grid torus by construction, see Figure \ref{fig:3x3gridtorus}.
	
	\begin{figure}[ht]
		\centering
		\includegraphics[width=0.3\textwidth]{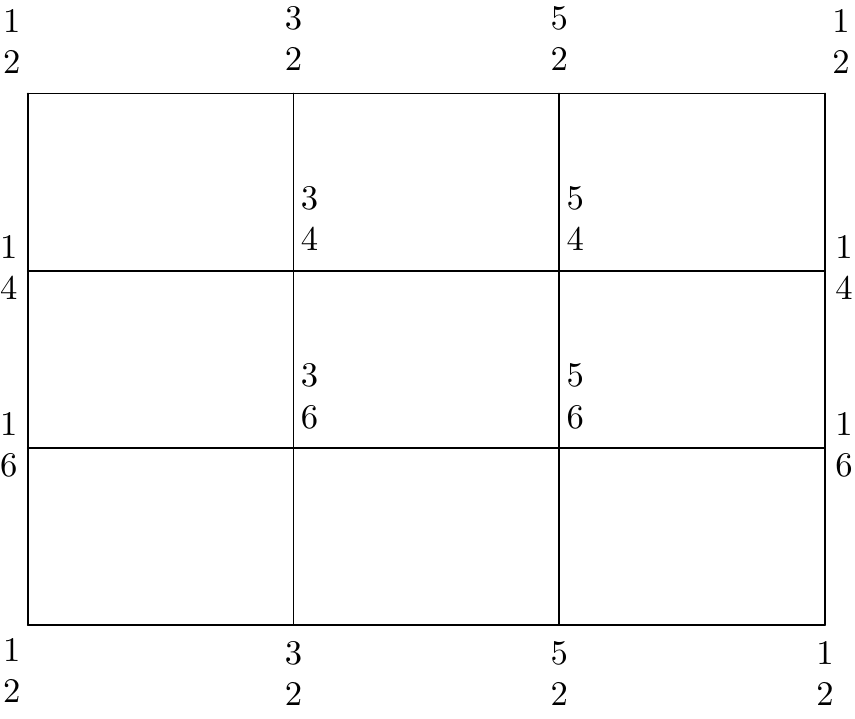}
		\caption{The standard $3 \times 3$-grid torus as a quadrangulated discrete normal surface of genus $1$ of the combinatorial $3$-sphere $ BdC_4(6) $.}
		\label{fig:3x3gridtorus}
	\end{figure}
\end{example}

Looking at slicings $S_{(V_1,V_2)}$ where either $\dim (\operatorname{span} (V_1) ) = 1$ or $\dim (\operatorname{span} (V_2) ) = 1$ we can prove the following for arbitrary combinatorial $3$-manifolds:

\begin{satz}
	\label{thm:quadrsli}
	Let $M$ be a combinatorial $3$-manifold, $S_{(V_1,V_2)}$ a connected slicing of genus $g$ of $M$, $q$ its number of quadrilaterals, $\dim(\operatorname{span} (V_i)) = 1$ and  $n:=|V_i|$ for $i = 1$ or $i = 2$. Then
	\begin{equation}
		\label{eq:quadrsli}
		q \geq 3(n+g-1).
	\end{equation}
	\noindent
	In particular, this applies to all quadrangulated slicings.
\end{satz}

\begin{proof}
	To obtain a connected slicing $S_{(V_1,V_2)}$ of genus $g$ both the span of $V_1$ and of $V_2$ must be connected and contain a graph of genus $g$. The Euler characteristic tells us that a graph of genus $g$ with $n$ vertices needs exactly $(g+n-1)$ edges ($n \geq k$ for $ \binom{k-2}{2} < g \leq \binom{k-1}{2}$). Note that every edge of the manifold is surrounded by at least three tetrahedra. Since $\dim (\operatorname{span} (V_1) ) = 1$, each of the at least $3(g+n-1)$ tetrahedra must be distinct from all the others. In particular, exactly one edge of each of these tetrahedra lies in $\operatorname{span} (V_1)$ and, thus, the intersection of $S_{(V_1,V_2)}$ with each of the tetrahedra appears as a quadrilateral in the slicing.
	
	If a slicing is quadrangulated, then $\dim (\operatorname{span} (V_1) ) = \dim (\operatorname{span} (V_2) ) = 1$ and the conditions for Theorem \ref{thm:quadrsli} are fulfilled.
\end{proof}

\begin{figure}[h]
	\centering
	\includegraphics[width=0.4\textwidth]{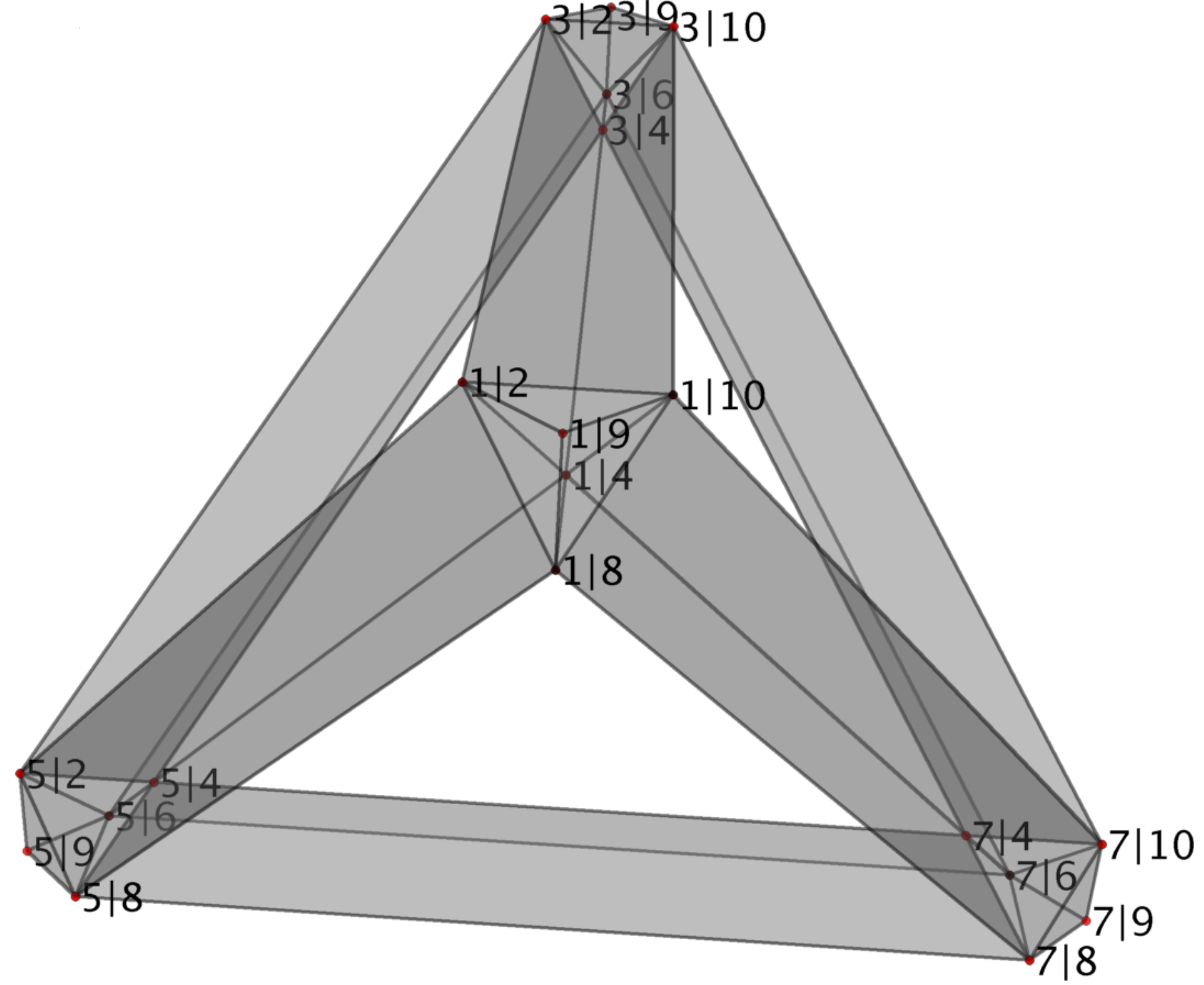}
	\includegraphics[width=0.5\textwidth]{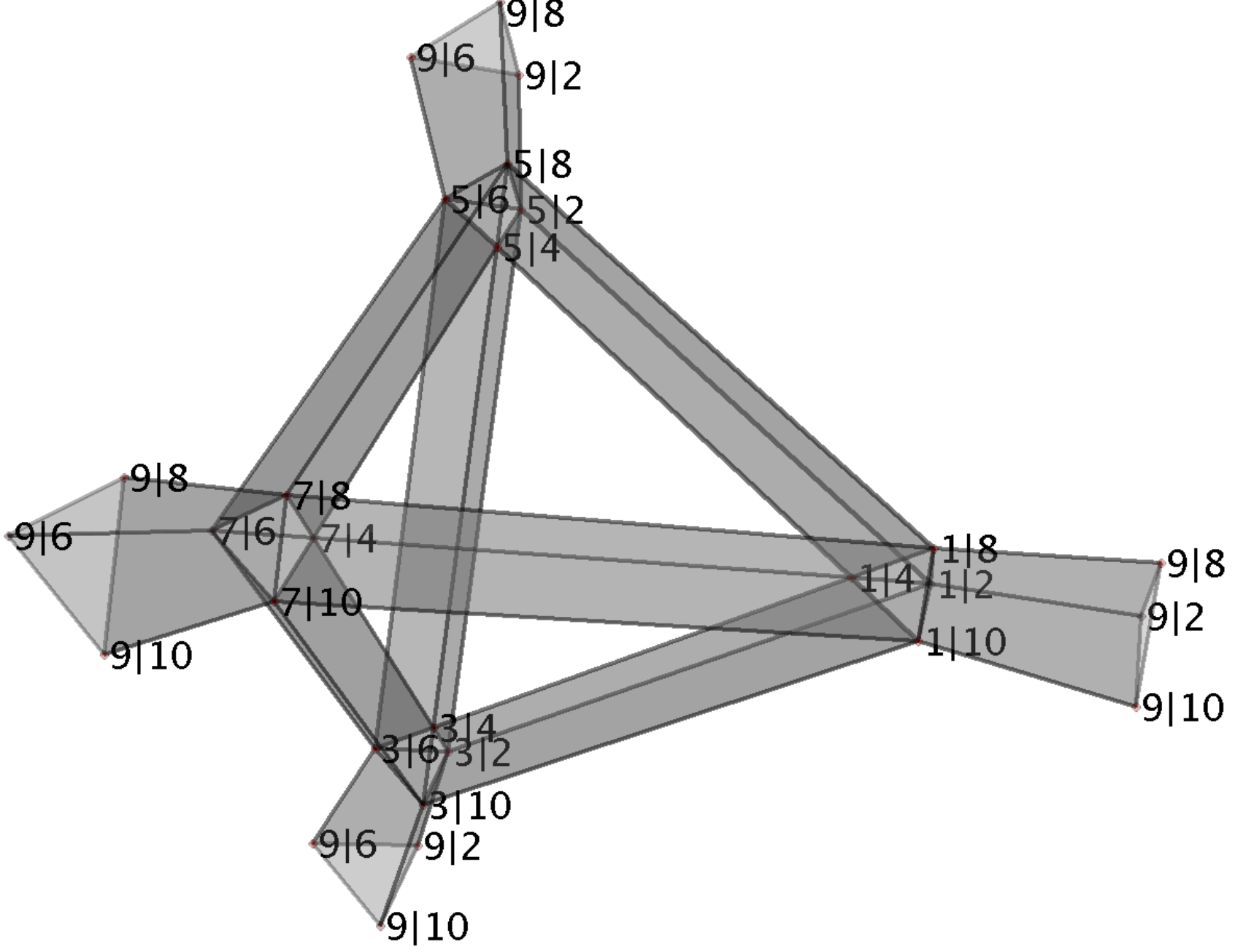}
	\caption{Slicings of the nearly neighborly centrally symmetric $10$-vertex
	sphere $ ^3 10^{22}_{1}$ (cf. \cite{Lutz03TrigMnfFewVertVertTrans}). On the left: Slicing behind the complete $1$-skeleton of $\langle 1,3,5,7 \rangle$ - a surface of genus $3$.	On the right: Slicing between the vertices $V_1 := \{ 1,3,5,7,9 \}$ and $V_2 := \{ 2,4,6,8,10 \}$. Both $V_1$ and $V_2$ span a complete graph $K_{5,5}$, generating a quadrangulated surface of genus $6$.}
	\label{fig:slicing4-5}
\end{figure}

For some examples of slicings that attain equality in (\ref{eq:quadrsli}) consider the nearly neighborly centrally symmetric $4$-polytope $P$ with $10$ vertices and the maximum number of $40$ edges and transitive automorphism group 
\begin{equation*}
	\operatorname{Aut} (P) = \langle (1,4,7,6,9,2)(3,10)(5,8), (1,7,3,9)(2,8,4,6) \rangle \cong C_2 \times S_5
\end{equation*}
due to Gr\"unbaum (cf. Section $6.4$ of \cite{Gruenbaum67ConvPoly}). The boundary $\partial P = S $ is a	combinatorial $3$-sphere ($ ^3 10^{22}_{1}$ in \cite{Lutz03TrigMnfFewVertVertTrans}) and equals the $ \operatorname{Aut} (P) $-orbit $$S = (1,2,3,4)_{30}$$ of length $30$. $S$ admits a rsl-function $f$ such that every slicing induced by $f$ attains equality in (\ref{eq:quadrsli}): 
\begin{center}
	\begin{tabular}{l|c|c|c}
		slicing & $ g $ & $n$ & $ \# $ quadrilaterals  \\
		\hline
		$\{ 1 \}, \{ 2,3,4,5,6,7,8,9,10\}$ & $ 0 $ & $1$ & $0$ \\
		$\{ 1,3 \}, \{ 2,4,5,6,7,8,9,10\}$ & $ 0 $ & $2$ & $3$ \\
		$\{ 1,3,5 \}, \{ 2,4,6,7,8,9,10 \}$ & $ 1 $ & $3$ & $9$ \\
		$\{ 1, 3, 5, 7 \},\{ 2, 4, 6, 8, 9, 10 \}$ & $4$ & $ 3 $ & $18$ \\
		$\{ 1, 3, 5, 7, 9 \},\{ 2, 4, 6, 8, 10 \}$ & $5$ & $ 6 $ & $30$ \\
	\end{tabular}
\end{center} 
For a picture of some of the slicings see Figure \ref{fig:slicing4-5}.

\bigskip
In addition, there is a $15$-vertex triangulation $C$ of $S^2 \dtimes S^1$ with transitive automorphism group 
\begin{equation*}
	\operatorname{Aut} (C) = \langle (1, \ldots , 15 ) \rangle \cong \mathbb{Z}_{15}
\end{equation*}
($ ^{3} 15^{1}_{17}$ in \cite{Lutz03TrigMnfFewVertVertTrans}), generated by the $\mathbb{Z}_{15}$-orbits	$$ C = (1,2,3,5)_{15} \cup (1,2,3,12)_{15} \cup (1,2,4,6)_{15} \cup (1,2,5,7)_{15} \cup (1,2,6,7)_{15}, $$ which admits a number of examples of non-orientable surfaces that are obtained via slicings that fulfill equality in (\ref{eq:quadrsli}):
\begin{center}
	\begin{tabular}{l|c|c|c}
		slicing & $ g $ & $n$ & $ \# $ quadrilaterals  \\
		\hline
		$\{ 1,4,7 \},$ & $ 1 $ & $ 3 $ & $ 9 $ \\ $ \{ 2,3,5,6,8, \ldots ,15 \}$&&&\\
		$\{ 1,4,7,10 \},$ & $ 3 $ & $ 4 $ & $ 18 $ \\ $ \{ 2,3,5,6,8,9,11,12,13,14,15 \}$&&&\\
		$\{ 1,4,7,10,13 \},$ & $ 6 $ & $ 5 $ & $ 30 $ \\$ \{ 2,3,5,6,8,9,11,12,14,15 \}$&&&\\
		\\
	\end{tabular}
\end{center} 

The examples were found by an enumerative search using the support for slicings  of the \texttt{GAP} package \texttt{simpcomp} \cite{simpcomp}. The complexes that were used to construct the slicings are included in the built-in library of \texttt{simpcomp}.

\begin{remark}
	\label{rem:conj}
	It seems natural to assume that (\ref{eq:quadrsli}) also holds in the generic case of arbitrary combinatorial $3$-manifolds since the condition $\dim (\operatorname{span} (V_1) ) = 1$ only holds for examples of low complexity. However, this is not true:
	
An analogue of	Theorem \ref{thm:quadrsli} in the general case would state, that at least $39$ quadrilaterals are needed for a surface of genus $8$. However, there exist orientable and non-orientable slicings of genus $8$ with only $35$ quadrilaterals (slicing $(\{ 1, 3, 5, 7, 9, 11, 13 \},\{ 2, 4, 6, 8, 10, 12, 14 \} )$ of $^3 14^1_{3}$ and $(\{ 1, 3, 5, 7, 9, 11, 13 \},\{ 2, 4, 6, 8, 10, 12, 14 \})$ of $^3 14^1_{32}$, notations as in \cite{Lutz03TrigMnfFewVertVertTrans}). Moreover, there is an orientable slicing of genus $2$ with only $14$ quadrilaterals whereas the lower bound in Theorem \ref{thm:quadrsli} would be $15$ ($(\{ 1, 2, 16, 21 \},\{ 3, 4, 5, 6, 7, 8, 9, 10, 11, 12,	13, 14, 15, 17, 18, 19, 20, 22 \})$ of a triangulation of $(S^2 \dtimes S^1)^{\# 14}$ from \cite{Lutz08ManifoldPage}).	
\end{remark}

As mentioned in Remark \ref{rem:conj} a straightforward generalization of Theorem \ref{thm:result} and Theorem \ref{thm:quadrsli} does not seem to be obvious. Nevertheless, we state:

\begin{conj}
	\label{conj:main}
	Let $M$ be a combinatorial $3$-manifold and $S$ a slicing of $M$ with $q$ quadrilaterals. Then 
	\begin{equation*}
		q \geq 3 - \frac32 \chi(S)
	\end{equation*}
	holds.
\end{conj}

Purely combinatorial methods do not seem to be suitable to prove Conjecture \ref{conj:main}. However, at least for low genera a geometric or algebraic approach could lead to new results. 

Furthermore it remains to investigate if the stated theorems or conjectures can be extended to discrete normal surfaces. Although most proves are founded on the fact that a slicing splits the surrounding combinatorial manifold into two pieces it is believed that Theorem \ref{thm:result} and Theorem \ref{thm:quadrsli} are also true for discrete normal surfaces that do not define such a splitting.

\section{Weakly neighborly slicings}
\label{sec:neighslic}

From the Dehn-Sommerville equations for $2$-manifolds it follows, that the $f$-vector of a triangulated surface $S$ is already determined by the number of vertices $n$ of $S$:
\begin{equation}
	f(S) = \left ( n, 3(n- \chi (S)), 2(n- \chi (S) \right ).
\end{equation}
Since obviously the number of edges cannot exceed ${ n \choose 2 }$, this gives rise to the well known lower bound on the number of vertices needed to triangulate $S$:
\begin{eqnarray}
	\label{eq:heawood}
	n &\geq& \frac12 \left ( 7 + \sqrt{49 - 24 \chi (S)} \right ).
\end{eqnarray}
Inequality (\ref{eq:heawood}) is often referred to as {\it Heawood inequality}. The inequality is sharp whenever $3(n- \chi(S)) = { n \choose 2 }$. Note that not all topological types of surfaces admit $2$-neighborly triangulations but whenever they do this triangulation is minimal with respect to the number of vertices (these are the so called {\it regular cases} of the Heawood inequality). In contrast to the simplicial case, the term of a neighborly complex does not make sense for a polyhedral map since there are pairs of vertices contained in one face which do not span an edge. This directly leads us to the following.

\begin{deff}[Weakly neighborly polyhedral map, from \cite{Brehm86WeaklyNeighbPolyMaps}]
	We call a polyhedral map {\it weakly neighborly} if for any two vertices there is a face containing both of them.
\end{deff}

Together with the above definition we can give a lower bound on the number of vertices for polyhedral $m$-gon maps:

\begin{prop}[cf. Lemma 3 in \cite{Brehm86WeaklyNeighbPolyMaps}]
	\label{thm:mgonwithoutdiags}
	Let $S$ be a polyhedral $m$-gon map with Euler characteristic $\chi ( S )$, then $S$ needs at least
	\begin{equation}
		\label{eq:mgonwithoutdiags}
		n \geq \frac12 \left ( (2m+1)  + \sqrt{(2m + 1)^2 - 8m \chi(S)} \right ) 
	\end{equation}
	vertices, with equality if and only if $S$ is weakly neighborly.	
\end{prop}

\begin{proof}
	Proposition \ref{thm:mgonwithoutdiags} is a special case of \cite[Lemma 3]{Brehm86WeaklyNeighbPolyMaps} with $n = V$ and $m = k_i$ for all $1 \leq i \leq p = F$, where $p = \frac{2}{m-2} (n - \chi (S) )$ denotes the number of $m$-gons of $S$.
%
\end{proof}

In the case $m = 4$ we have
\begin{equation}
	\label{eq:quadwithoutdiags}
	n \geq \frac12 \left ( 9  + \sqrt{81 - 32 \chi(S)} \right ). 
\end{equation}
Thus, the $9$-vertex grid torus from Figure \ref{fig:3x3gridtorus} is an example for a weakly neighborly quadrangulation of the torus. 

Let us now consider the case of polyhedral maps consisting of triangles and quadrilaterals.

\begin{lemma}[cf. Lemma 2 in \cite{Brehm86WeaklyNeighbPolyMaps}]
	\label{lemma:equiv}
	Let $S$ be a polyhedral map with $n$ vertices only consisting of $t$ triangles and $q$ quadrilaterals, then the following statements are equivalent:
	\begin{enumerate}
		\item[{\it a)}] $S$ is weakly neighborly,
		\item[{\it b)}] $S$ has $ e = { n \choose 2 } - 2 q $ edges,
		\item[{\it c)}] $n = \frac12 \left (7 + \sqrt{49 + 8  q - 24 \chi (S) } \right )$,
		\item[{\it d)}]	$q = { n-3 \choose 2 } +  3 \chi (S) - 6$.
	\end{enumerate}
\end{lemma}

\begin{proof}
	\begin{description}
		\item[$ a) \Leftrightarrow b) $] Every quadrilateral in $S$ decreases the	maximal number of edges in $S$ by $2$. In particular, $S$ has exactly $ e = { n 	\choose 2 } - 2 q $ edges. This can be seen as a special case of  \cite[Lemma 2 (5)]{Brehm86WeaklyNeighbPolyMaps} with $n = V$, $t$ times $k_i = 3$ and $q$ times  $k_i = 4$. 
		\item[$ b) \Leftrightarrow c) $] By the pseudomanifold property we have $ 2e = 3 t + 4 q $ and with the Euler characteristic $ \chi (S)  = n - e + \frac{2e - 4 q}{3} + q $	and $ e =  { n \choose 2 } - 2 q$ we get
		\begin{equation*}
			\chi (S) = n - \frac13 \left( { n \choose 2 } - 3 q \right ),
		\end{equation*}
		which is equivalent to
		\begin{equation*}
			0 = n^2 - 7n - 2 q + 6 \chi (S).
		\end{equation*}
		It follows that
		\begin{equation}
			\label{eq:genHeawood}
			n = \frac12 \left (7 \pm \sqrt{49 + 8 q - 24 \chi (S)} \right ) ,
		\end{equation}
		\noindent
		where clearly only the greater one of the two solutions of the quadratic equation is valid (see also \cite[Lemma 2 (6)]{Brehm86WeaklyNeighbPolyMaps} for a similar statement).
		\item[$ c) \Leftrightarrow d) $] Given $\chi (S)$ and $n$ we have at most one non-negative integer $q$ solving (\ref{eq:genHeawood}). Thus, we get 
		\small
		\begin{eqnarray}
			n &=& \frac12 (2n) \nonumber \\
			 &=& \frac12 \left(7 + \sqrt{(2n-7)^2} \right ) \nonumber \\
			 &=& \frac12 \left(7 + \sqrt{49 + 4(n^2 -7n) + 24 \chi (S) - 24 \chi (S) } \right )
			\nonumber \\
			 &=& \frac12 \left(7 + \sqrt{49 + 8 \left (\frac{n^2 -7n}{2} +6-6 + 3 \chi (S) \right ) - 24 \chi (S) } \right ) \nonumber \\
			 &=& \frac12 \left(7 + \sqrt{49 + 8 \left ({ n-3 \choose 2} - 6 + 3 \chi (S) \right ) - 24 \chi (S) } \right ). \nonumber
		\end{eqnarray}
		\normalsize
		This directly leads to
		\begin{equation*}
			q = { n-3 \choose 2} - 6 + 3 \chi (S).
		\end{equation*}
	\end{description}
\end{proof}

In a series of papers \cite{Brehm85WeaklPolyhMapsTorus,Brehm86NonexWNPMapsGenus2,Brehm86WeaklyNeighbPolyMaps}
Altshuler and Brehm classified weakly neighborly maps on surfaces.

\medskip
It is obvious (for example by looking at the examples of  \cite{Brehm85WeaklPolyhMapsTorus}), that not all weakly neighborly polyhedral maps can be realized as discrete normal surfaces or slicings of combinatorial pseudomanifolds. In fact, this case is a rare exception. In order to see that, just note that the number of possible edges of a discrete normal surface is strongly restricted by Lemma \ref{lemma:neighbors}. More precisely we have the following condition:

\begin{lemma}
	Let $S_{(V_1,V_2)}$ be a weakly neighborly slicing of a combinatorial $3$-pseudomanifold $M$ induced by the vertex partition $V_1 \dot{\cup} V_2 = V$, $n_i$ the number of vertices that lie in the boundary of $\operatorname{span}(V_i)$, $i=1,2$, then
	\begin{equation}
		\label{eq:maxneighb0}
		n_1 n_2 ( 15 - n_1 n_2 - n_1 - n_2 ) = 12 \chi (S_{(V_1,V_2)}).
	\end{equation}
\end{lemma}

\begin{proof}
	We first need to verify that any pair of vertices $a \in V_1$, $y \in V_2$ in the boundary of $\operatorname{span} (V_1)$, $\operatorname{span} (V_2)$ occurs  as a vertex ${ a \choose y }$ of $S_{(V_1,V_2)}$: Since $a$ ($y$) lies in the boundary of $\operatorname{span} (V_1)$ ($\operatorname{span} (V_2)$) there exist two vertices $b \in V_1$ and $x \in V_2$ such that $\langle a,x \rangle, \, \langle b,y \rangle \in M$. If $a=b$ or $x=y$ then ${ a \choose y }$ is a vertex of $S_{(V_1,V_2)}$. Now let $a \neq b$ and $x\neq y$ then ${a \choose x}$ and ${b \choose y}$ are both vertices of $S_{(V_1,V_2)}$. By Lemma \ref{lemma:neighbors} they can not form an edge of $S_{(V_1,V_2)}$. However, since $S_{(V_1,V_2)}$ is weakly neighborly, ${a \choose x}$ and ${b \choose y}$ lie in one facet and, thus, $\langle { a \choose x}, {a \choose y}, { b \choose x}, {b \choose y} \rangle $ must be one of its quadrialterals. As a consequence ${a \choose y}$ is a vertex of $S_{(V_1,V_2)}$ and the slicing has $n_1 n_2$ vertices.	
	
	No it remains to show that the boundary of $\operatorname{span} (V_i)$, $ i = 1,2$, is $2$-neighborly:	Let $a,b \in V_1$ be two vertices in the boundary of $\operatorname{span}(V_1)$, $a \neq b$. Then there exist two edges $ \langle a , x \rangle$ and $ \langle b , y \rangle$ with $x,y \in V_2$. If $x \neq y$, then $\langle { a \choose x}, {b \choose y} \rangle$ can not be an edge of $S_{(V_1,V_2)}$ by Lemma \ref{lemma:neighbors}. Since $S_{(V_1,V_2)}$ is weakly neighborly $\langle { a \choose x}, {a \choose y}, { b \choose x}, {b \choose y} \rangle $ must be a quadrilateral, $\langle a,b,x,y \rangle$ a tetrahedron of $M$ and $\langle a,b \rangle$ an edge of $\operatorname{span}(V_1)$. If $x=y$, then $\langle { a \choose x}, {b \choose x} \rangle$ must be an edge of $S_{(V_1,V_2)}$ since it is weakly neighborly. Hence, $\{ a,b,x\}$ spans a triangle in $M$ and $\langle a,b \rangle$ is an edge. Altogether, the boundary of $\operatorname{span}(V_i)$, $i \in \{ 1,2 \}$ must be $2$-neighborly and $S_{(V_1,V_2)}$ must have exactly 
	\begin{equation}
		\label{eq:maxneighb}
		n_1 { n_2 \choose 2 } + n_2 { n_1 \choose 2 } = { n_1 n_2 \choose 2 } - 2 \left( { n_1 n_2 - 3 \choose 2 } + 3 \chi (S) - 6 \right)
	\end{equation}
	edges, where the right hand side follows from Lemma \ref{lemma:equiv} $d)$.	This yields	(\ref{eq:maxneighb0}).
\end{proof}

If we furthermore restrict $M$ to be a combinatorial manifold we can classify all such weakly neighborly slicings:

\begin{satz}
	The only weakly neighborly polyhedral maps that are slicings of a combinatorial $3$-manifold $M$ are the boundary of the $3$-simplex, the boundary of a triangular prism with $6$ vertices and the $3 \times 3$-grid torus shown in Figure \ref{fig:3x3gridtorus}.	
\end{satz}

\begin{proof}
	Since $M$ is a $3$-manifold with $\chi (M) = 0$, $\chi (S)$ must be an even number. 
	
	In the case $\chi (S) = 2$ the only integer solutions of (\ref{eq:maxneighb0}) with $\chi (S) = 2$ are $n_1 = 1$, $n_2 = 4$ and $n_1 =2$, $n_2 = 3$. In the first case we have a (triangulated) vertex figure which has to be $2$-neighborly. Hence, the only triangulated $2$-neighborly $2$-sphere is the boundary of the $3$-simplex. The second case must occur as a slicing of a closed combinatorial manifold with $2+3=5$ vertices (since the span of $3$ vertices can not have interior faces in dimension $3$). The only combinatorial $3$-manifold with $5$ vertices is the boundary of the $4$-simplex. The slicing between an edge and a triangle yields a weakly neighborly $2$-sphere consisting of $3$ quadrialterals and two triangles: this has to be the boundary of a $3$-dimensional prism $\Delta^1 \times \Delta^2$.
	
	In the case $\chi (S) = 0$: Since $M$ is a combinatorial $3$-manifold, it follows from Proposition \ref{prop:handles} that $n_1 , n_2 \geq 3 $. The only solution is $n_1 = n_2 = 3$. Since the span of $3$ vertices can not have interior faces in dimension $3$, the slicing must have $6$ vertices (which necessarily forms a triangulated $3$-sphere by virtue of Theorem A in \cite{Brehm87CombMnfFewVert}) and two disjoint empty triangles. This determines the combinatorial type of the complex and of the slicing. The unique solution is the $3 \times 3$-grid torus shown in Figure \ref{fig:3x3gridtorus}. In particular, the Klein bottle can not appear as a weakly neighborly slicing (by the orientability of $S^3$).
	
	In the case $\chi (S) \leq -2$ the lower bounds on $n_1$ and $n_2$ given in Proposition \ref{prop:handles} and the asymptotic behaviour of the left hand side of (\ref{eq:maxneighb0}) do not admit further solutions.	
\end{proof}

\addcontentsline{toc}{chapter}{Bibliography}

{\footnotesize
 \bibliographystyle{abbrv}
 \bibliography{/home/spreerjn/LaTeX/dissbib/bibliography.bib}

\begin{thebibliography}{10}

\bibitem{Brehm86NonexWNPMapsGenus2}
A.~Altshuler and U.~Brehm.
\newblock Nonexistence of weakly neighborly polyhedral maps on the orientable
  {$2$}-manifold of genus {$2$}.
\newblock {\em J. Combin. Theory Ser. A}, 42(1):87--103, 1986.

\bibitem{Barnette71MinNumVertSimplePoly}
D.~Barnette.
\newblock {T}he minimum number of vertices of a simple polytope.
\newblock {\em Israel J. Math.}, 10:121--125, 1971.

\bibitem{Barnette73ProofLBCConvPoly}
D.~Barnette.
\newblock {A} proof of the lower bound conjecture for convex polytopes.
\newblock {\em Pacific J. Math.}, 46:349--354, 1973.

\bibitem{Brehm85WeaklPolyhMapsTorus}
U.~Brehm and A.~Altshuler.
\newblock The weakly neighborly polyhedral maps on the torus.
\newblock {\em Geom. Dedicata}, 18(3):227--238, 1985.

\bibitem{Brehm86WeaklyNeighbPolyMaps}
U.~Brehm and A.~Altshuler.
\newblock On weakly neighborly polyhedral maps of arbitrary genus.
\newblock {\em Israel J. Math.}, 53(2):137--157, 1986.

\bibitem{Brehm87CombMnfFewVert}
U.~Brehm and W.~K{\"u}hnel.
\newblock {C}ombinatorial manifolds with few vertices.
\newblock {\em Topology}, 26(4):465--473, 1987.

\bibitem{Brehm93PolyhedralMflds}
U.~Brehm and J.~M. Wills.
\newblock Polyhedral manifolds.
\newblock In P.~M. Gruber and J.~M. Wills, editors, {\em Handbook of convex
  geometry, {V}ol.\ {A}, {B}}, pages 535--554. North-Holland, Amsterdam, 1993.

\bibitem{simpcomp}
F.~Effenberger and J.~Spreer.
\newblock {simpcomp - A GAP package}, {V}ersion 1.4.
\newblock \url{http://www.igt.uni-stuttgart.de/LstDiffgeo/simpcomp}, 2010.
\newblock {S}ubmitted to the {\it GAP Group}.

\bibitem{Gruenbaum67ConvPoly}
B.~Gr{\"u}nbaum.
\newblock {\em {C}onvex polytopes}.
\newblock With the cooperation of Victor Klee, M. A. Perles and G. C. Shephard.
  Pure and Applied Mathematics, Vol. 16. Interscience Publishers John Wiley \&
  Sons, Inc., New York, 1967.

\bibitem{Haken61TheoNormFl}
W.~Haken.
\newblock {T}heorie der {N}ormalfl{\"a}chen.
\newblock {\em Acta Math.}, 105:245--375, 1961.

\bibitem{Haken62HomeomProb3Mflds}
W.~Haken.
\newblock \"{U}ber das {H}om\"oomorphieproblem der $3$-{M}annigfaltigkeiten.
  {I}.
\newblock {\em Math. Z.}, 80:89--120, 1962.

\bibitem{Kalai87RigidityLBT}
G.~Kalai.
\newblock {R}igidity and the lower bound theorem. {I}.
\newblock {\em Invent. Math.}, 88(1):125--151, 1987.

\bibitem{Kalelkar08ChiQuadrNS}
T.~Kalelkar.
\newblock {E}uler characteristic and quadrilaterals of normal surfaces.
\newblock {\em Proc. Indian Acad. Sci. Math. Sci.}, 118(2):227--233, 2008.

\bibitem{Kneser29ClosedSurfIn3Mflds}
H.~Kneser.
\newblock {Geschlossene Fl\"{a}chen in dreidimensionalen Mannigfaltigkeiten}.
\newblock {\em Jahresbericht der deutschen Mathematiker-Vereinigung},
  38:248--260, 1929.

\bibitem{Kuhnel86MinTrigKummVar}
W.~K{\"u}hnel.
\newblock {M}inimal triangulations of {K}ummer varieties.
\newblock {\em Abh. Math. Sem. Univ. Hamburg}, 57:7--20, 1986.

\bibitem{Kuehnel90TrigMnfFewVert}
W.~K{\"u}hnel.
\newblock {T}riangulations of manifolds with few vertices.
\newblock In {\em Advances in differential geometry and topology}, pages
  59--114. World Sci. Publ., Teaneck, NJ, 1990.

\bibitem{Kuehnel95TightPolySubm}
W.~K{\"u}hnel.
\newblock {\em {T}ight polyhedral submanifolds and tight triangulations},
  volume 1612 of {\em Lecture Notes in Mathematics}.
\newblock Springer-Verlag, Berlin, 1995.

\bibitem{Kuehnel84RhombiTess3Space15Vertex3Torus}
W.~K{\"u}hnel and G.~Lassmann.
\newblock {T}he rhombidodecahedral tessellation of {$3$}-space and a particular
  {$15$}-vertex triangulation of the {$3$}-dimensional torus.
\newblock {\em Manuscripta Math.}, 49(1):61--77, 1984.

\bibitem{Kuehnel85NeighbComb3MfldsDihedralAutGroup}
W.~K{\"u}hnel and G.~Lassmann.
\newblock {N}eighborly combinatorial {$3$}-manifolds with dihedral automorphism
  group.
\newblock {\em Israel J. Math.}, 52(1-2):147--166, 1985.

\bibitem{Lickorish}
W.~B.~R. Lickorish.
\newblock {\em An introduction to knot theory}, volume 175 of {\em Graduate
  Texts in Mathematics}.
\newblock Springer-Verlag, New York, 1997.
\newblock Kapitel 12.

\bibitem{Lutz08ManifoldPage}
F.~H. Lutz.
\newblock {T}he {M}anifold {P}age.
\newblock {\url{http://www.math.tu-berlin.de/diskregeom/stellar}}.

\bibitem{Lutz03TrigMnfFewVertVertTrans}
F.~H. Lutz.
\newblock {\em {T}riangulated {M}anifolds with {F}ew {V}ertices and
  {V}ertex-{T}ransitive {G}roup {A}ctions}.
\newblock {S}haker {V}erlag, Aachen, 1999.
\newblock PhD Thesis, TU Berlin.

\bibitem{Spreer09CombPorpsOfK3}
J.~Spreer and W.~K{\"u}hnel.
\newblock {C}ombinatorial properties of the {K}3 surface: {S}implicial blowups
  and slicings.
\newblock {\tt arXiv:0909.1453v2 [math.CO]}, preprint, 31 pages, 3 figures,
  2009.
\newblock To appear in Exp. Math.

\bibitem{Tillmann07NormalSurfacesTopFinite3Mflds}
S.~Tillmann.
\newblock Normal surfaces in topologically finite 3-manifolds.
\newblock {\em Enseign. Math. (2)}, 54(3-4):329--380, 2008.

\end{thebibliography}
}

\noindent
Institut f\"ur Geometrie und Topologie \\
Universit\"at Stuttgart \\
70550 Stuttgart \\
Germany

\end{document}